\documentclass[a4paper,11pt]{scrartcl}
\usepackage{amsmath,amsthm,amssymb}
\usepackage{enumerate}
\usepackage{graphicx}
\usepackage[square,sort&compress,numbers]{natbib}
\usepackage{mathtools}\mathtoolsset{showonlyrefs=false}
\usepackage{subcaption}
\newcommand{\Authornote}{\renewcommand{\thefootnote}{\fnsymbol{footnote}}}
\newcommand{\authornote}{\Authornote\footnote}

\theoremstyle{plain}
\newtheorem{theorem}{Theorem}

\theoremstyle{definition}

\newtheorem{algorithm}{Algorithm}
\newtheorem{example}{Example}
\theoremstyle{remark}
\newtheorem{remark}[theorem]{Remark}

\newcommand{\refex}[1]{Example~\ref{#1}}

\newcommand{\refalg}[1]{Algorithm~\ref{#1}}

\newcommand{\refrem}[1]{Remark~\ref{#1}}
\newcommand{\reffig}[1]{Figure~\ref{#1}}

\newcommand{\finbox}{\nolinebreak\hfill{\small $\blacksquare$}}

\newcounter{alnum}
\newenvironment{algstep}
{\begin{list}{{\upshape\bfseries Step~\arabic{alnum}:}}%
{\usecounter{alnum}%
\setlength{\leftmargin}{18mm}%
\setlength{\labelwidth}{20mm}%
}}{\end{list}}

\newcommand{\Min}{\mathop{\mathrm{Min.}}}
\newcommand{\Max}{\mathop{\mathrm{Max.}}}
\newcommand{\st}{\mathop{\mathrm{s.{\,}t.}}}

\newcommand{\kernel}{\mathop{\mathrm{Ker}}\nolimits}
\newcommand{\diag}{\mathop{\mathrm{diag}}\nolimits}

\renewcommand{\Re}{\ensuremath{\mathbb{R}}}

\newcommand{\bi}[1]{\ensuremath{\boldsymbol{#1}}}

\newcommand{\rr}[1]{\ensuremath{\mathrm{#1}}}

\begin{document}

\begin{center}
  {\Large\bfseries\sffamily% 
  Redundancy Optimization of Finite-Dimensional }
  \par\medskip
  {\Large\bfseries\sffamily% 
  Structures: A Concept and a Derivative-Free Algorithm 
  }%
  \par\bigskip
  Yoshihiro Kanno \authornote[2]{
    Laboratory for Future Interdisciplinary Research of Science and Technology, 
    Institute of Innovative Research, 
    Tokyo Institute of Technology, 
    Nagatsuta 4259, Yokohama 226-8503, Japan.
    E-mail: \texttt{kanno.y.af@m.titech.ac.jp}. 
    Phone: +81-45-924-5364. 
    Fax: +81-45-924-5978.
  }
\end{center}

\begin{abstract}
  Redundancy is related to the amount of functionality that the structure 
  can sustain in the worst-case scenario of structural degradation. 
  This paper proposes a widely-applicable concept of redundancy 
  optimization of finite-dimensional structures. 
  The concept is consistent with the robust structural optimization, as 
  well as the quantitative measure of structural redundancy based on the 
  information-gap theory. 
  A derivative-free algorithm is proposed based on the sequential 
  quadratic programming (SQP) method, where we use the finite-difference 
  method with adaptively varying the difference increment. 
  Preliminary numerical experiments show that an optimal solution of 
  the redundancy optimization problem possibly has multiple worst-case 
  scenarios. 
\end{abstract}
\begin{quote}
  \textbf{Keywords}
  \par
  Robustness; 
  redundancy; 
  resilience; 
  uncertainty; 
  derivative-free optimization; 
  simplex gradient. 
\end{quote}

\section{Introduction}

Since real-world structures inevitably encounter various uncertainties, 
robustness and redundancy are crucial concepts in structural design. 
Robust optimization of structures has been studied extensively; see, e.g., 
\cite{BtN97,CC03,GZZ13,HK15,JvDvK11,Sig09,TNKK11,YK10}, and the 
references therein. 
In contrast, study on design methodology of structures considering 
redundancy property is very limited. 
\citet{MSMK14} applied the notion of redundancy in the coding theory to 
truss optimization. 
\citet{MG14} considered the conditional probability of failure of each 
structural component given failure of the structure, and proposed to 
minimize the difference between the maximum and minimum values 
of these conditional probabilities. 
Also, optimization problems of fail-safe structures are relevant to 
the problem studied in this paper. 
\citet{SAH76} and \citet{NA82} defined a damage condition as 
complete or partial damage to selected structural components. 
They specified several damage scenarios a priori, and 
imposed the performance constraints on all the damaged scenarios, as 
well as on the intact scenario. 
For continuum-based topology optimization, \citet{JLSS14} considered 
the effect of local failure. 
As a simple model of local failure, they supposed that the material 
stiffness in a patch with a specified shape vanishes. 
The location of patch is assumed to be unknown in advance, and the 
compliance in the worst-case scenario is minimized to obtain 
a fail-safe structure. 

This paper attempts to present a widely-applicable concept of redundancy 
optimization of finite-dimensional structures. 
Several definitions of redundancy of structures have been proposed; see, 
e.g., \cite{BB99,BL13,FC87,KBh11,KHKS16,SB05,ZF12,Zih00}, 
and the references therein. 
Among others, in this paper we adopt the strong redundancy \citep{KBh11}, 
which is a quantitative measure of structural redundancy based on the 
information-gap theory \citep{Bh06}. 
The strong redundancy of a given structural design is defined as the 
greatest deficiency that can arise at any place in the structure, without 
violating the performance constraint. 
Thus the concept of strong redundancy can fully address uncertainty in 
future structural degradation. 
Moreover, it has flexibility to incorporate diverse structural 
performances, unlike most of other redundancy measures that 
consider only the ultimate strength or the collapse load as 
the structural performance of interest. 

In the redundancy optimization problem proposed in this paper, we 
specify only the upper bound for the number of damaged structural 
components, denoted $\alpha$, while the damaged structural components 
themselves are considered unknown (or, uncertain). 
Then we formulate a worst-case design optimization problem. 
Namely, we optimize the objective function evaluated in the worst-case 
scenario of structural degradation. 
This methodology shares a common framework with extensively studied 
robust compliance optimization of structures with a non-probabilistic 
modeling of uncertainty. 
There, the worst-case compliance (i.e., the maximum value of the 
compliance) is minimized, when the static external load or the 
structural geometry is assumed to be uncertain  
\citep{BtN97,CC03,GZZ13,HK15,JvDvK11,Sig09,TNKK11,YK10}. 

As a concrete example of the redundancy optimization, attention of this 
paper is focused on maximization of the worst-case limit load factor 
of a truss structure. 
When $\alpha=0$, i.e., no redundancy is required, the redundancy 
optimization reverts to the classical limit design 
(optimal plastic design). 
For $\alpha \ge 1$, an optimal solution of the redundancy optimization 
is statically indeterminate, unlike the limit design. 
For given $\alpha$ and structural design, the worst-case limit load 
factor can be computed via mixed-integer linear programming 
(MILP) \citep{Kan12}; see also \citet{TTlWG13}. 
As a substitute of the gradient of the worst-case limit load factor, 
we employ the simplex gradient 
(also called the stencil gradient) that is often used in derivative-free 
optimization methods \citep{CSV09,Kel11}. 
More concretely, we compute the finite-difference gradient with 
decreasing a difference increment as the optimization procedure 
progresses. 
Making use of this approximated gradient, we propose a derivative-free 
method based on the sequential quadratic programming (SQP) 
for solving the redundancy optimization problem.

The paper is organized as follows. 
First, a concept of redundancy optimization of structures is defined. 
We also discuss the relation between the redundancy optimization and the 
robust optimization. 
The following section  presents a derivative-free optimization 
method that combines SQP and the finite-difference gradient with a 
varying difference increment. 
Then, we apply the algorithm to simple problem instances to 
investigate some properties of the obtained solutions. 
Some conclusions are drawn at the end of the paper.

%\paragraph{Notation}

A few words regarding notation. 
We use ${}^{\top}$ to denote the transpose of a vector or a matrix. 
For two vectors $\bi{x} =(x_{i}) \in \Re^{n}$ and 
$\bi{s} =(s_{i}) \in \Re^{n}$, 
we write $\bi{x} \ge \bi{s}$ if $x_{i} \ge s_{i}$ $(i=1,\dots,n)$. 
Particularly,  $\bi{x} \ge \bi{0}$ means $x_{i} \ge 0$ $(i=1,\dots,n)$. 
We use $\bi{1} = (1,1,\dots,1)^{\top}$ to denote the all-ones vector. 
The $\ell_{1}$-norm and the Euclidean norm of a vector 
$\bi{x} \in \Re^{n}$ is denoted by 
$\| \bi{x} \|_{1} = \sum_{i=1}^{n} |x_{i}|$ and 
$\| \bi{x} \| = (\bi{x}^{\top} \bi{x})^{1/2}$, respectively. 
We use $\diag(\bi{x})$ to denote the $n \times n$ diagonal matrix with a 
vector $\bi{x} \in \Re^{n}$ on its diagonal. 
For a matrix $A \in \Re^{m \times n}$, 
we use $A^{+} \in \Re^{n \times m}$ and $\kernel A$ to denote 
its Moore--Penrose pseudoinverse and nullspace, respectively. 
For $a$, $b \in \Re$ satisfying $a < b$, we denote by $[a,b]$ and 
$]a,b[$ the closed and open intervals between $a$ and $b$, respectively.

\section{Problem formulation}

In this section we present a concept of redundancy optimization of 
structures, which is consistent with the widely accepted notion of 
robust structural optimization \citep{BtN97,CC03,HK15,TNKK11,YK10}. 
Also, it is naturally endowed with a qualitative measure of structural 
redundancy, called the strong redundancy \citep{KBh11}. 

Let $\bi{x} \in \Re^{m}$ denote the vector of design variables, where 
$m$ is the number of the design variables.  
Consider structural performance $h(\bi{x})$ depending on structural 
design $\bi{x}$. 
A small value of $h(\bi{x})$ is preferred over a large value. 
The conventional optimization problem that maximizes the structural 
performance is formulated as follows: 
\begin{subequations}\label{P.nominal.optimization}
  \begin{alignat}{3}
    & \Min
    &{\quad}& 
    h(\bi{x}) \\
    & \st && 
    \bi{x} \in X . 
  \end{alignat}
\end{subequations}
Here, $X \subseteq \Re^{m}$ is the set of feasible design variables. 

\begin{example}\label{ex.truss.1}
  Consider a truss structure consisting of $m$ members. 
  Let $\bi{x}$ be a nonnegative real vector, the $i$th component of 
  which is the cross-sectional area of member $i$. 
  We use $\lambda(\bi{x})$ to denote the limit load factor of the truss, 
  and let $h(\bi{x})=-\lambda(\bi{x})$. 
  A small value of $h(\bi{x})$ is worthy. 
  A typical example of constraint $\bi{x} \in X$ has the form 
  \begin{align}
    X = \{ \bi{x} 
    \mid
    \bi{c}^{\top} \bi{x} \le V, 
    \
    \bi{x} \ge \bi{0} 
    \} , 
    \label{eq.truss.example.X}
  \end{align}
  where $c_{i}$ is the undeformed length of member $i$ and $V$ is the 
  specified upper bound for the structural volume. 
  In this situation, problem \eqref{P.nominal.optimization} is the 
  conventional limit design problem. 
  \finbox
\end{example}

\begin{example}\label{ex.displacement}
  As in \refex{ex.truss.1}, consider a truss structure. 
  Let $\bi{u}(\bi{x})$ denote the nodal displacement vector caused by 
  the specified static external load. 
  Consider the displacement constraints 
  \begin{align}
    |u_{j}(\bi{x})| \le \bar{g}_{j} , 
    \quad  j=1,\dots,k , 
    \label{eq.ex.displacement.constraint.1}
  \end{align}
  where $\bar{g}_{j}$ is a specified positive value and $k$ is the 
  number of nodes for which the displacement constraint is considered. 
  Let $h$ be 
  \begin{align*}
    h(\bi{x}) 
    = \max \{ |u_{1}(\bi{x})| - \bar{g}_{1} , \dots, 
    |u_{k}(\bi{x})| - \bar{g}_{k} , 0 \} . 
  \end{align*}
  Namely, $h(\bi{x})$ measures violation of constraint 
  \eqref{eq.ex.displacement.constraint.1}, and $h(\bi{x})=0$ means that 
  constraint \eqref{eq.ex.displacement.constraint.1} is satisfied. 
  \finbox
\end{example}

Future structural damage is unknown in advance. 
Hence, based upon the information-gap theory \citep{Bh06}, the strong 
redundancy \citep{KBh11} of design $\bi{x}$ is defined as the greatest 
deficiency that can arise at any place in the structure, without 
violating the performance constraint. 
In other words, the strong redundancy guarantees structural 
functionality in the worst-case scenario, when future structural damage 
is uncertain. 
Therefore, to increase redundancy it is natural to attempt to maximize 
the structural performance in the worst-case scenario of structural 
deficiency. 
Obviously, the worst-case scenario depends on the structural design. 

For structural component $i$, we use binary variable $t_{i}$ that serves 
as a indicator of soundness. 
Specifically, the value of $t_{i}$ is defined by 
\begin{align*}
  t_{i} = 
  \begin{dcases*}
    1 & if member $i$ is intact, \\
    0 & if member $i$ is damaged. \\
  \end{dcases*}
\end{align*}
Then vector $\bi{t} = (t_{1},\dots,t_{m})^{\top}$ expresses the 
scenario of deficiency that the structure suffers. 
In particular, the nominal scenario, which refers to the completely 
intact structure, corresponds to $\bi{t}=\bi{1}$. 
Following \citet{Kan12} and \citet{KBh11}, we define the deficiency set by 
\begin{align*}
  T(\alpha) 
  = \{ \bi{t} \in \{ 0,1 \}^{m}
  \mid 
  \| \bi{t} - \bi{1} \|_{1} \le \alpha
  \} , 
\end{align*}
where $\alpha \ge 0$ is a parameter representing the level of structural 
damage. 
Namely, $T(\alpha)$ is the set of all scenarios in which the structure 
suffers degradation of at most an amount $\alpha$. 
From the definition, it is straightforward to see 
that $T(0)$ is a singleton consisting of the nominal scenario (i.e., 
$\bi{t}=\bi{1}$), and that $0 \le \alpha \le \alpha'$ implies 
$T(\alpha) \subseteq T(\alpha')$. 

Recall that the structural design is characterized by $\bi{x}$. 
We assume that a damaged structural component completely loses its 
functionality. 
Then the realization of the $i$th structural component is written as 
$t_{i}x_{i}$. 
Therefore, the set of all possible realizations of the structural design 
is given by 
\begin{align}
  D(\bi{x}; \alpha) 
  = \{ \diag(\bi{t}) \bi{x} 
  \mid 
  \bi{t} \in T(\alpha)
  \} . 
  \label{eq.def.deficiency.set.1}
\end{align}

\begin{remark}\label{rem.partial.deficiency}
  In \eqref{eq.def.deficiency.set.1} we assume a model that a damaged 
  structural component is completely missing from the structure. 
  Alternatively, we may suppose that structural components are 
  diminishing only in part. 
  Let $\gamma \in [0,1[$ be a constant representing the degree of damage. 
  We assume that all members share same value of $\gamma$. 
  Then the deficiency set is given by 
  \begin{align*}
    D(\bi{x}; \alpha) 
    = \{ \diag(\bi{t} + \gamma (\bi{1}-\bi{t})) \bi{x} 
    \mid 
    \bi{t} \in T(\alpha)
    \} . 
  \end{align*}
  This model with $\alpha=1$ was considered in, e.g., 
  \cite{SAH76}. 
  Obviously, when $\gamma=0$, this model reverts to 
  \eqref{eq.def.deficiency.set.1}. 
  \finbox
\end{remark}

For given $\alpha \ge 0$ and $\bi{x} \in X$, define 
$h^{\rr{worst}}(\bi{x}; \alpha)$ by 
\begin{align}
  h^{\rr{worst}}(\bi{x}; \alpha) 
  = \max \{ h(\bi{s}) \mid \bi{s} \in D(\bi{x}; \alpha) \}  . 
  \label{eq.def.h.worst}
\end{align}
Namely, $h^{\rr{worst}}(\bi{x}; \alpha)$ is the value of structural 
performance when the structure suffers the worst-case damage scenario. 
When the amount of structural degradation, $\alpha \ge 0$, is specified, 
we attempt to improve the performance in the worst-case scenario as far 
as possible. 
This design optimization problem is formulated as follows: 
\begin{subequations}\label{P.redundancy.optimization}
  \begin{alignat}{3}
    & \Min
    &{\quad}& 
    h^{\rr{worst}}(\bi{x}; \alpha) \\
    & \st && 
    \bi{x} \in X . 
  \end{alignat}
\end{subequations}
It is worth noting that, when $\alpha=0$, problem 
\eqref{P.redundancy.optimization} reverts to problem 
\eqref{P.nominal.optimization}, i.e., the conventional optimization 
problem. 
The amount of uncertainty increases as $\alpha$ increases. 
In the following, we call problem \eqref{P.redundancy.optimization} the 
redundancy optimization problem. 

Problem \eqref{P.redundancy.optimization} is maximization of the 
objective function evaluated at the worst-case scenario, when the set of 
damaged members is unknown. 
This can be viewed as a robust optimization problem; see \citet{BtEgN09} 
for the notion of robust optimization. 
For instance, robust compliance optimization of structures, that has been 
studied extensively, attempts to minimize the compliance at the 
worst-case scenario, when the external load 
and/or the structural geometry are not known precisely
\citep{BtN97,CC03,TNKK11,YK10,GZZ13,HK15,JvDvK11,Sig09}.

\section{Derivative-free SQP method}

In this section, we develop an algorithm for solving the redundancy 
counterpart of the concrete problem discussed in \refex{ex.truss.1}. 
That is, $\bi{x}$ is the vector of member cross-sectional areas, 
$-h(\bi{x})$ is the limit load factor, and $X$ is defined by 
\eqref{eq.truss.example.X}. 
However, the algorithm presented below may be applicable to a broader 
class of problems in the form \eqref{P.redundancy.optimization}. 

The limit load factor of a truss is defined as follows. 
Let $d$ denote the number of degrees of freedom of the nodal 
displacements. 
Suppose that the external load consists of a constant part, denoted 
$\bi{p}_{\rr{d}}$, and a proportionally increasing part, denoted 
$\lambda\bi{p}_{\rr{r}}$, where $\bi{p}_{\rr{d}} \in \Re^{d}$ and 
$\bi{p}_{\rr{r}} \in \Re^{d}$ are constant vectors and 
$\lambda \in \Re$ is a load factor. 
We use $\bi{b}_{i} \in \Re^{d}$ $(i=1,\dots,m)$ to denote the $i$th 
column vector of the equilibrium matrix. 
It follows from the lower bound theorem of the limit analysis that the 
limit load factor is the optimal value of the following linear 
programming problem: 
\begin{subequations}\label{P.lower.limit.analysis}
  \begin{alignat}{3}
    & \Max
    &{\quad}& 
    \lambda \\
    & \st && 
    \sum_{i=1}^{m} q_{i} \bi{b}_{i} 
    = \lambda \bi{p}_{\rr{r}} + \bi{p}_{\rr{d}} , 
    \label{P.lower.limit.analysis.2} \\
    & && 
    |q_{i}| \le \sigma_{\rr{y}} x_{i} , 
    \quad  i=1,\dots,m. 
  \end{alignat}
\end{subequations}
Here, $\lambda$ and $q_{1},\dots,q_{m}$ are variables to be optimized, 
and $\sigma_{\rr{y}}>0$ is the (constant) yield stress. 

For given $\alpha \ge 0$ and $\bi{x} \ge \bi{0}$, it is known that the 
value of $h^{\rr{worst}}(\bi{x}; \alpha)$ defined by 
\eqref{eq.def.h.worst} can be computed via MILP \citep{Kan12}. 
However, $h^{\rr{worst}}(\bi{x}; \alpha)$ is not necessarily 
differentiable with respect to $\bi{x}$. 
As an approximation of the gradient, if any, of 
$h^{\rr{worst}}(\bi{x}; \alpha)$, or as its substitute, 
we employ the stencil gradient, that is often used in 
derivative-free optimization methods \citep{CSV09,Kel11}, 
as an approximation of the gradient of the objective function. 
By making use of the stencil gradient, we construct 
a quadratic programming problem that approximates the original problem 
in \eqref{P.redundancy.optimization}. 

Let 
\begin{align*}
  f(\bi{x}) := h^{\rr{worst}}(\bi{x}; \alpha) 
\end{align*}
for notational simplicity. 
Then the problem to be solved is written as follows: 
\begin{subequations}\label{P.to.be.solved}
  \begin{alignat}{3}
    & \Min
    &{\quad}& 
    f(\bi{x}) \\
    & \st && 
    \bi{c}^{\top} \bi{x} \le V , 
    \label{P.to.be.solved.1} \\
    & && 
    \bi{x} \ge \bi{0} . 
  \end{alignat}
\end{subequations}

We begin with computation of the stencil gradient 
(also called a simplex gradient) of $f$, denoted $\nabla_{\rr{s}}f$. 
Essentially we approximate the gradient by using the finite difference 
method and, as is done in the implicit filtering method, reduce 
the difference increment as the optimization procedure progresses; 
see, e.g., \citet{CSV09} and \citet{Kel11} for fundamentals of the 
stencil gradient and the implicit filtering method. 
Since the implicit filtering uses a relatively large difference 
increment at the early stage of optimization, it may possibly neglect 
the high-frequency low-amplitude features of the objective function and 
avoid the algorithm converging to a poor local optimal 
solution \citep{Kel11}. 
From \eqref{P.lower.limit.analysis}, we can see that, 
if $\bi{x} \le \bi{x}'$, then the limit load factor of $\bi{x}'$ is no 
less than that of $\bi{x}$. 
Therefore, constraint \eqref{P.to.be.solved.1} becomes active at an 
optimal solution. 
This motivates us to use only points satisfying \eqref{P.to.be.solved.1} 
with equality as sample points for the finite-difference method. 
Let $\{ \bi{\delta}_{1}, \dots, \bi{\delta}_{m-1} \}$ 
be a basis of $\kernel \bi{c}^{\top}$, 
where $\bi{\delta}_{i} \in \Re^{m}$ is normalized as 
$\| \bi{\delta}_{i} \| = 1$ $(i=1,\dots,m-1)$. 
Define the set of sample points, centered at $\bi{x}$, by 
\begin{align}
  S(\bi{x};r) 
  = \{ \bi{x} \pm r \bi{\delta}_{i} 
  \mid i=1,\dots,m-1 \} , 
  \label{eq.sample.points.set.1}
\end{align}
where constant $r > 0$ is called the stencil radius. 
If a sample point defined above has a negative component, then we modify 
it so as to be a feasible point; 
see \refrem{rem.sample.point} for more accounts. 
For notational simplicity, we use $\bi{z}_{i}$ to denote an element of 
$S(\bi{x};r)$, i.e., 
\begin{align*}
  S(\bi{x};r) = \{ \bi{z}_{1}, \dots, \bi{z}_{2m-2} \} . 
\end{align*}
Define $Y \in \Re^{(2m-2) \times m}$ and $\bi{\delta} \in \Re^{2m-2}$ by 
\begin{align*}
  Y &= 
  \begin{bmatrix}
    \bi{z}_{1}-\bi{x}, \dots, \bi{z}_{2m-2}-\bi{x} 
  \end{bmatrix}
  ^{\top} , \\
  \bi{\delta} &= 
  \begin{bmatrix}
    f(\bi{z}_{1}) - f(\bi{x}) \\
    \vdots \\
    f(\bi{z}_{2m-2}) - f(\bi{x}) \\
  \end{bmatrix}
  . 
\end{align*}
Then the stencil gradient of $f$ at the point $\bi{x}$ is given 
as follows \citep{CSV09,Kel11}: 
\begin{align}
  \nabla_{\rr{s}} f(\bi{x}) = Y^{+} \bi{\delta} . 
  \label{eq.def.simplex.gradient}
\end{align}

\begin{remark}\label{rem.sample.point}
  A sample point defined by \eqref{eq.sample.points.set.1} may possibly 
  have a negative member cross-sectional area. 
  At such a sample point, the limit load factor is not defined and, 
  hence, the value of $f$ is not defined. 
  Therefore, we replace negative cross-sectional areas with a small 
  constant $\varepsilon > 0$ and reduce positive cross-sectional areas so 
  that the structural volume of the resulting sample point becomes 
  equal to $V$. 
  Consequently, all elements of $S(\bi{x};r)$ are positive vectors 
  satisfying \eqref{P.to.be.solved.1} with equality. 
  \finbox
\end{remark}

Making use of $\nabla_{\rr{s}}f(\bi{x})$ in 
\eqref{eq.def.simplex.gradient}, we next design a derivative-free SQP 
method for solving problem \eqref{P.to.be.solved}. 
Let $\bi{x}_{k}$ denote the incumbent solution at the $k$th iteration. 
We solve the following quadratic programming (QP) problem 
in variables $\bi{d} \in \Re^{m}$ to determine the search direction: 
\begin{subequations}\label{eq.P.SQP.subproblem.1}
  \begin{alignat}{3}
    & \Min
    &{\quad}& 
    \frac{1}{2} \bi{d}^{\top} B_{k} \bi{d} 
    + \nabla_{\rr{s}} f(\bi{x}_{k})^{\top} \bi{d} \\
    & \st && 
    \bi{c}^{\top} \bi{d} \le V - \bi{c}^{\top} \bi{x}_{k}, \\
    & && 
    \bi{d} \ge -\bi{x}_{k} . 
  \end{alignat}
\end{subequations}
Here, $B_{k}$ is a symmetric positive definite matrix. 
Let $\bi{d}_{k}$ denote an optimal solution of problem 
\eqref{eq.P.SQP.subproblem.1}. 
We employ $\bi{d}_{k}$ as a search direction, and perform the line 
search to determine the step length, denoted $a_{k}$. 
Then the incumbent solution, $\bi{x}_{k}$, is updated as 
\begin{align*}
  \bi{x}_{k+1} = \bi{x}_{k} + a_{k} \bi{d}_{k} . 
\end{align*}

As matrix $B_{k}$ in \eqref{eq.P.SQP.subproblem.1}, we adopt a 
quasi-Newton approximation of the Hessian of the Lagrangian of 
problem \eqref{P.to.be.solved}. 
Here, the Lagrangian is defined by 
\begin{align*}
  L(\bi{x}, \mu, \bi{\zeta}) = 
  f(\bi{x}) 
  - \mu (V - \bi{c}^{\top} \bi{x}) 
  - \bi{\zeta}^{\top} \bi{x} , 
\end{align*}
where $\mu \ge 0$ and $\bi{\zeta} \ge \bi{0}$ are the Lagrange 
multipliers. 
More concretely, we employ the damped BFGS update formula, which is one 
of conventional formulae used in SQP \cite[section~18.3]{NW06}, to generate 
$B_{k+1}$ from $B_{k}$. 
Namely, we first compute $\bi{s}_{k}$, $\bi{y}_{k} \in \Re^{m}$ defined 
by 
\begin{align}
  \bi{s}_{k} &= \bi{x}_{k+1} - \bi{x}_{k} , \\
  \bi{y}_{k} &= \nabla_{\bi{x}}L(\bi{x}_{k+1},\mu_{k+1},\bi{\zeta}_{k+1}) 
  - \nabla_{\bi{x}}L(\bi{x}_{k},\mu_{k},\bi{\zeta}_{k}) . 
  \label{eq.BFGS.vector.y}
\end{align}
In practice, the gradient of the Lagrangian in \eqref{eq.BFGS.vector.y} 
is approximated by 
\begin{align*}
  \nabla_{\bi{x}}L(\bi{x}_{k}, \mu_{k}, \bi{\zeta}_{k}) 
  \simeq \nabla_{\rr{s}} f(\bi{x}_{k}) + \mu_{k} \bi{c} - \bi{\zeta}_{k} . 
\end{align*}
Also, $\mu_{k}$ and $\bi{\zeta}_{k}$ are approximated by the Lagrange 
multipliers of problem \eqref{eq.P.SQP.subproblem.1}. 
Next we compute $\theta_{k} \in \Re$ and $\bi{r}_{k} \in \Re^{m}$ by 
\begin{align*}
  \theta_{k} &= 
  \begin{dcases*}
    1 
    & if $\bi{y}_{k}^{\top} \bi{s}_{k} 
    \ge 0.2 \bi{s}_{k}^{\top} B_{k} \bi{s}_{k}$, \\
    0.8 \frac{\bi{s}_{k}^{\top} B_{k} \bi{s}_{k}}
    {\bi{s}_{k}^{\top} B_{k} \bi{s}_{k} - \bi{y}_{k}^{\top} \bi{s}_{k}}
    & otherwise, \\
  \end{dcases*}
  \\
  \bi{r}_{k} 
  &= \theta_{k} \bi{y}_{k} + (1-\theta_{k}) B_{k} \bi{s}_{k} . 
\end{align*}
Then $B_{k}$ is updated as follows: 
\begin{align}
  B_{k+1} 
  = B_{k} 
  - \frac{(B_{k} \bi{s}_{k})(B_{k} \bi{s}_{k})^{\top}}
  {\bi{s}_{k}^{\top} B_{k} \bi{s}_{k}}
  + \frac{\bi{r}_{k} \bi{r}_{k}^{\top}}{\bi{s}_{k}^{\top} \bi{s}_{k}} . 
  \label{eq.BFGS.formula}
\end{align}

We are now in position to describe the algorithm for solving 
problem \eqref{P.to.be.solved}. 

\begin{algorithm}[derivative-free SQP]
  \label{alg.SQP}
  {\quad}
  \begin{algstep}
    \setcounter{alnum}{-1}
    \item \label{alg.step.filter.0}
    Choose a feasible starting point $\bi{x}_{0}$. 
    Choose $r_{\min} > 0$, $r > r_{\min}$, $\tau_{\max} > 0$, 
    $\beta \in ]0,1[$, $\eta \in ]0,1[$, 
    $\rho \in ]0,1[$, $\epsilon > 0$, 
    and a symmetric positive definite matrix $B_{0} \in \Re^{m \times m}$. 
    Set $k := 0$. 
    
  \item \label{alg.step.filter.outer.2}
    If $r < r_{\min}$, then terminate. 

  \item \label{alg.step.filter.1}
    Generate a set of sample points, $S(\bi{x}_{k};r)$. 
    Compute the stencil gradient $\nabla_{\rr{s}}f(\bi{x}_{k})$ 
    by using the elements of $S(\bi{x}_{k};r)$. 
    If 
    \begin{align}
      \min \{ f(\bi{z}) \mid \bi{z} \in S(\bi{x}_{k};r) \}
       \ge f(\bi{x}_{k}) , 
      \label{eq.stencil.failure.1}
    \end{align}
    then set $r := \rho r$ and go to step~\ref{alg.step.filter.outer.2}. 

  \item \label{alg.SQP.subproblem}
    Solve problem \eqref{eq.P.SQP.subproblem.1}, and let $\bi{d}_{k}$ 
    denote the optimal solution. 
    If $\| \bi{d}_{k} \| < \epsilon$, then terminate. 
        
  \item \label{alg.SQP.line-search}
    Try to find the smallest integer $\tau \in [0,\tau_{\max}]$ 
    satisfying 
    \begin{align*}
      f(\bi{x}_{k} + \beta^{\tau} \bi{d}_{k}) 
      \le f(\bi{x}_{k})
      + \eta \beta^{\tau} \nabla_{\rr{s}}f(\bi{x}_{k})^{\top} \bi{d}_{k} .
    \end{align*}
    If such $\tau$ is successfully found, then let 
    $a_{k} := \beta^{\tau}$. 
    Otherwise, let $B_{k} := B_{0}$ and 
    $r := \rho r$, and go to step~\ref{alg.step.filter.outer.2}. 

  \item \label{alg.SQP.update.of.solution}
    Update $\bi{x}_{k}$ by 
    $\bi{x}_{k+1} := \bi{x}_{k} + a_{k} \bi{d}_{k}$. 
    Update $B_{k}$ to $B_{k+1}$ by \eqref{eq.BFGS.formula}. 
    Let $k \gets k+1$, and go to step~\ref{alg.step.filter.1}. 

  %\item 
    %If $\| \nabla_{\rr{s}}f(\bi{x}_{k})\| \le \tau r$, then 
    %set $r := r/2$ and go to step~\ref{alg.step.filter.outer.2}. 
    %If $k > k_{\max}$, then let $r := 2r$ and go to
    %step~\ref{alg.step.filter.outer.2}. 
    %Otherwise, let $k \gets k+1$, and go to 
    %step~\ref{alg.step.filter.1}. 
  \end{algstep}
\end{algorithm}

\begin{remark}
  At step~\ref{alg.step.filter.0} of \refalg{alg.SQP}, 
  we choose an initial point $\bi{x}_{0}$ 
  satisfying the constraints of problem \eqref{P.to.be.solved}. 
  On the other hand, we determine the step length $a_{k}$ 
  at step~\ref{alg.SQP.line-search} by performing 
  a conventional backtracking line search with the initial value $1$. 
  Moreover, the point $\bi{x}_{k}+\bi{d}_{k}$ is feasible for 
  problem \eqref{P.to.be.solved}. 
  Therefore, $\bi{x}_{k+1}$ determined at 
  step~\ref{alg.SQP.update.of.solution} is feasible for 
  problem \eqref{P.to.be.solved}. 
  %\memo{このことから，直線探索をするときに，メリット関数（の罰金項）は考えなくてよい．}
  %
  \finbox
\end{remark}

\begin{remark}
  At step~\ref{alg.step.filter.1}, we control the stencil radius, $r$, 
  according to the essential of the implicit filtering method 
  \citep{CSV09,Kel11}. 
  Namely, we decrease $r$ if the value of $f$ at the current point 
  $\bi{x}_{k}$ is no greater than the value at any sample point. 
  \finbox
\end{remark}

\begin{remark}\label{rem.optimality}
  At step~\ref{alg.SQP.subproblem}, we check a termination criterion. 
  If $f$ is differentiable at $\bi{x}_{k}$ and 
  $\nabla_{\rr{s}} f(\bi{x}_{k}) = \nabla f(\bi{x}_{k})$ holds, 
  then it follows from the fundamentals of the conventional SQP that 
  $\bi{d}_{k}=\bi{0}$ is a necessary condition for the local optimality 
  of problem \eqref{P.to.be.solved}. 
  Also, we might expect that $\nabla_{\rr{s}} f(\bi{x}_{k})$ 
  closely approximates $\nabla f(\bi{x}_{k})$ if $r$ is 
  sufficiently small. 
  The solutions found in the numerical experiments have multiple 
  worst-case scenarios. 
  At such a solution, the objective function, 
  $f(\bi{x})=h^{\rr{worst}}(\bi{x};\alpha)$ may not be differentiable. 
  The rigorous optimality condition for the redundancy optimization 
  remains to be studied. 
  \finbox
\end{remark}

\begin{remark}
  Step~\ref{alg.SQP.line-search} determines the step length in accordance 
  with the Armijo condition. 
  It is observed in the numerical experiments that $a_{k}=1$ is accepted 
  at many iterations. 
  However, since we use the stencil gradient instead of the gradient, 
  $\bi{d}_{k}$ computed at step~\ref{alg.SQP.subproblem} is not 
  necessarily a descent direction of $f$ at $\bi{x}_{k}$. 
  Therefore, it is possible that the line search at 
  step~\ref{alg.SQP.line-search} fails. 
  If this is the case, we decrease the stencil radius and compute again 
  the stencil gradient with a set of new sample points. 
  \finbox
\end{remark}

%\begin{itemize}
%  \item 
%\memo{``Stencil gradient'' という用語について\cite[p.~33]{CSV09}: 
%We note that simplex gradient when $p > n$ ($p$ is the number of 
%sample points) are also referred to as stencil gradients. 
%The set $\{ \bi{y}_{1},\dots,\bi{y}_{p} \}$ is a stencil centered 
%at $\bi{x}^{(k)}$. 
%For instance, the stencil could take the form 
%$\{ \bi{x}^{(k)} \pm h \bi{e}_{i} \mid i=1,\dots,n \}$, 
%where $h$ is the stencil radius : . 
%}
%\end{itemize}

\section{Preliminary numerical experiments}

In this section we apply \refalg{alg.SQP} to problem \eqref{P.to.be.solved}. 
We use simple problem instances to study some fundamental properties of 
the solutions obtained by the algorithm. 

\begin{figure}[tp]
  \centering
  \begin{subfigure}[b]{0.55\textwidth}
    \centering
    \scalebox{0.80}{
    \includegraphics{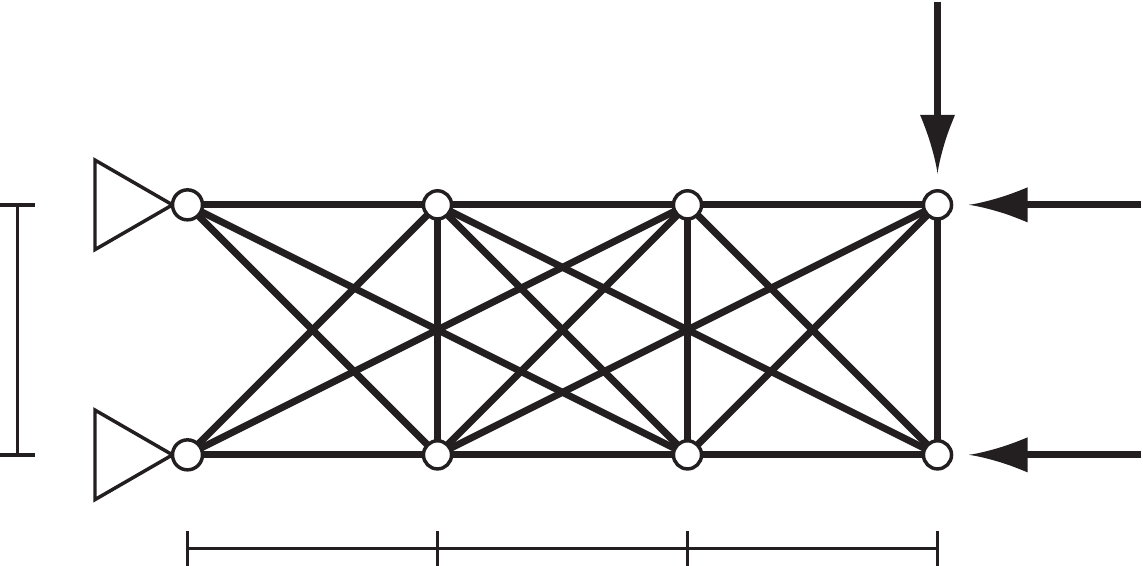}
    \begin{picture}(0,0)
      \put(-288,20){
      \put(-52,45){{\large $L$}}
      \put(228,120){{\large $\lambda p_{\rr{r}}$}}
      \put(260,92){{\large $p_{\rr{d}}$}}
      \put(260,2){{\large $p_{\rr{d}}$}}
      \put(42,-28){{\large $L$}}
      \put(114,-28){{\large $L$}}
      \put(186,-28){{\large $L$}}
      }
    \end{picture}
    }
    \caption{}
    \label{fig.m_19bar}
  \end{subfigure}
  \par\bigskip\bigskip
  \begin{subfigure}[b]{0.55\textwidth}
    \centering
    %%% C:\doc\derivative_free\redundancy\eva7\opt_design.m
    \includegraphics[scale=0.50]{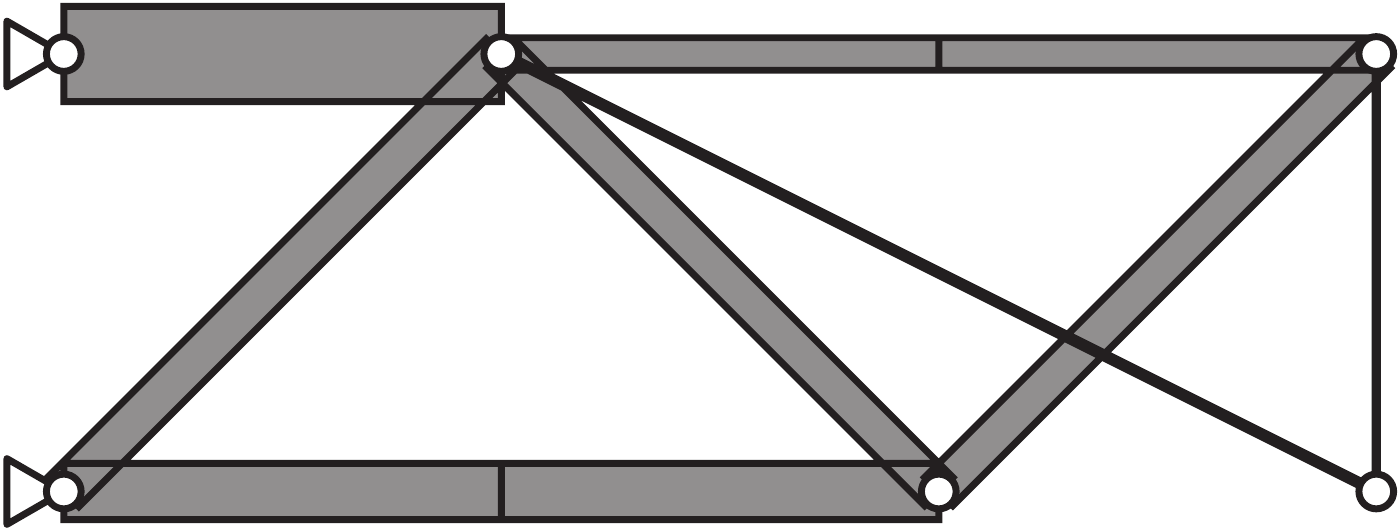}
    \caption{}
    \label{fig.optimal_alpha_0}
  \end{subfigure}
  \caption{Example (I). 
  \subref{fig.m_19bar}~Problem setting; and 
  \subref{fig.optimal_alpha_0}~the optimal solution without considering 
  redundancy (i.e., $\alpha=0$). 
  }
\end{figure}
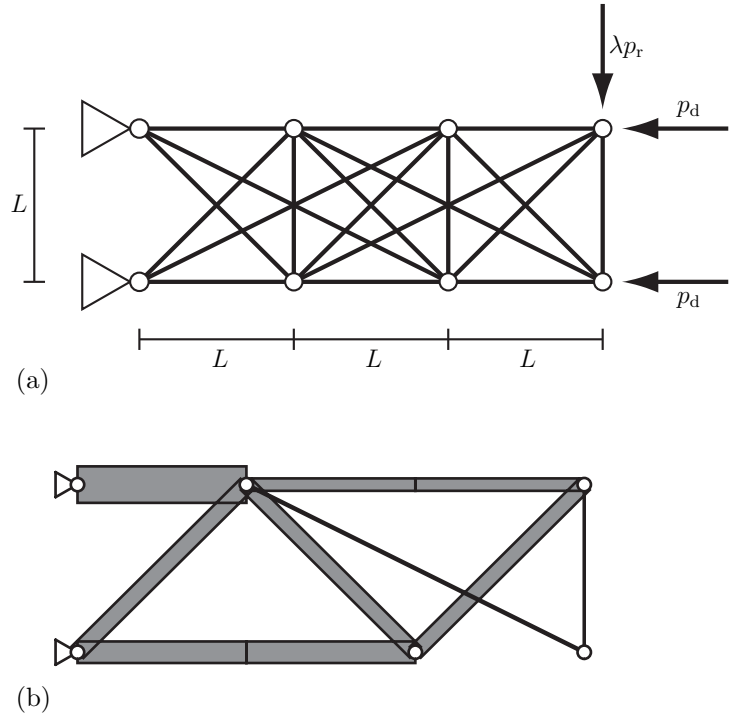

\begin{figure}[tp]
  \centering
  \begin{subfigure}[b]{0.55\textwidth}
    \centering
    %%% C:\doc\derivative_free\redundancy\eva6\opt_design.m
    \includegraphics[scale=0.50]{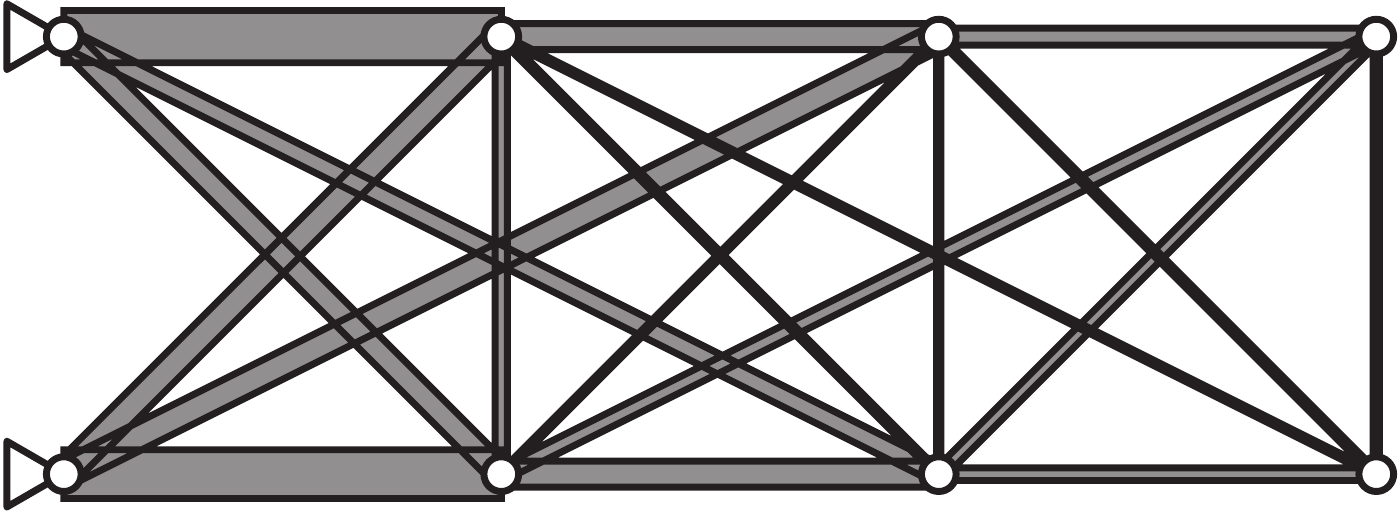}
    \caption{}
    \label{fig.optimal_alpha_1}
  \end{subfigure}
  \par\bigskip\bigskip
  \begin{subfigure}[b]{0.55\textwidth}
    \centering
    %%% C:\doc\derivative_free\redundancy\eva5\opt_design.m
    \includegraphics[scale=0.50]{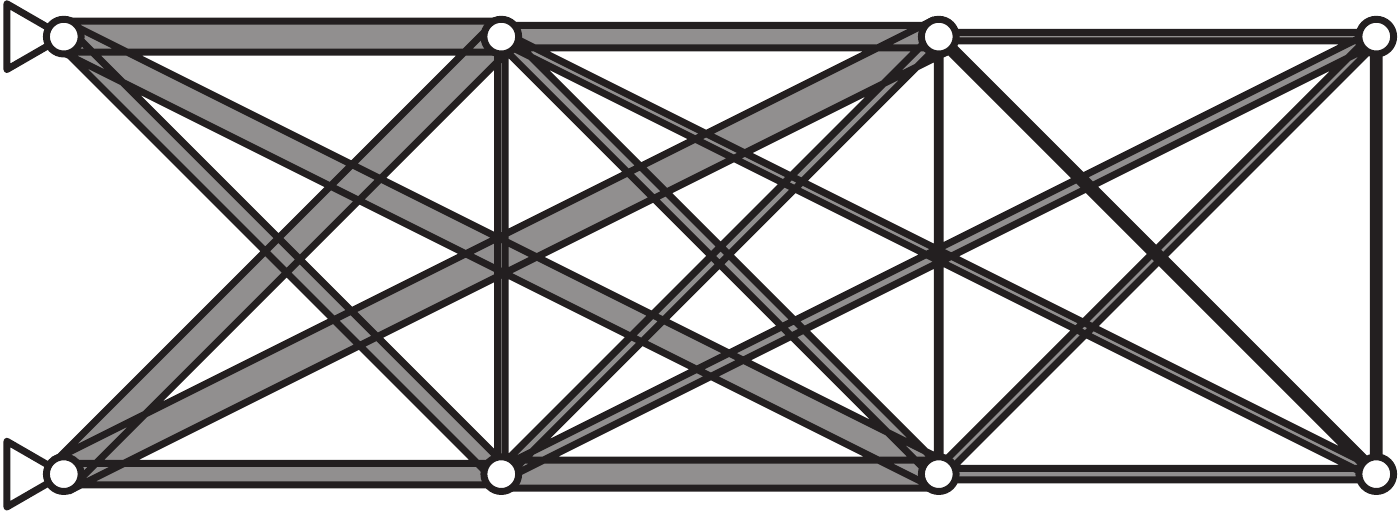}
    \caption{}
    \label{fig.optimal_alpha_2}
  \end{subfigure}
  \caption{The solutions obtained in example (I). 
  \subref{fig.optimal_alpha_1}~$\alpha=1$; and 
  \subref{fig.optimal_alpha_2}~$\alpha=2$.}
\end{figure}

\begin{figure}[tp]
  %%% C:\doc\derivative_free\redundancy\eva6\opt_design.m
  \centering
  \begin{subfigure}[b]{0.45\textwidth}
    \centering
    \includegraphics[scale=0.40]{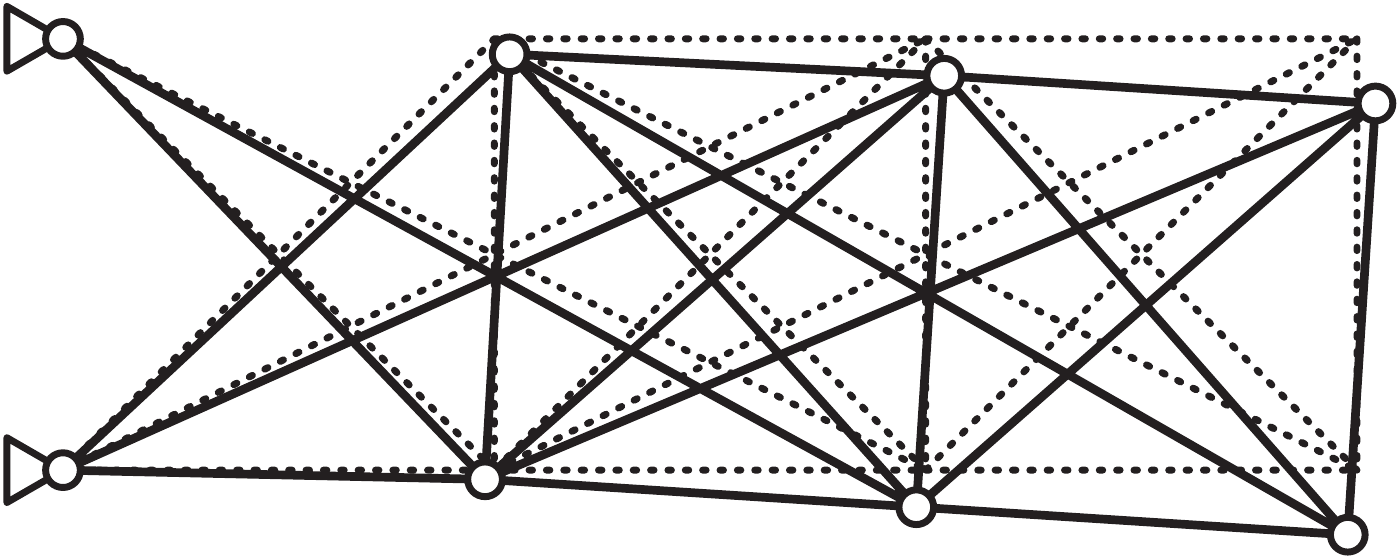}
    \caption{}
    \label{fig.enum1_alpha1}
  \end{subfigure}
  \hfill
  \begin{subfigure}[b]{0.45\textwidth}
    \centering
    \includegraphics[scale=0.40]{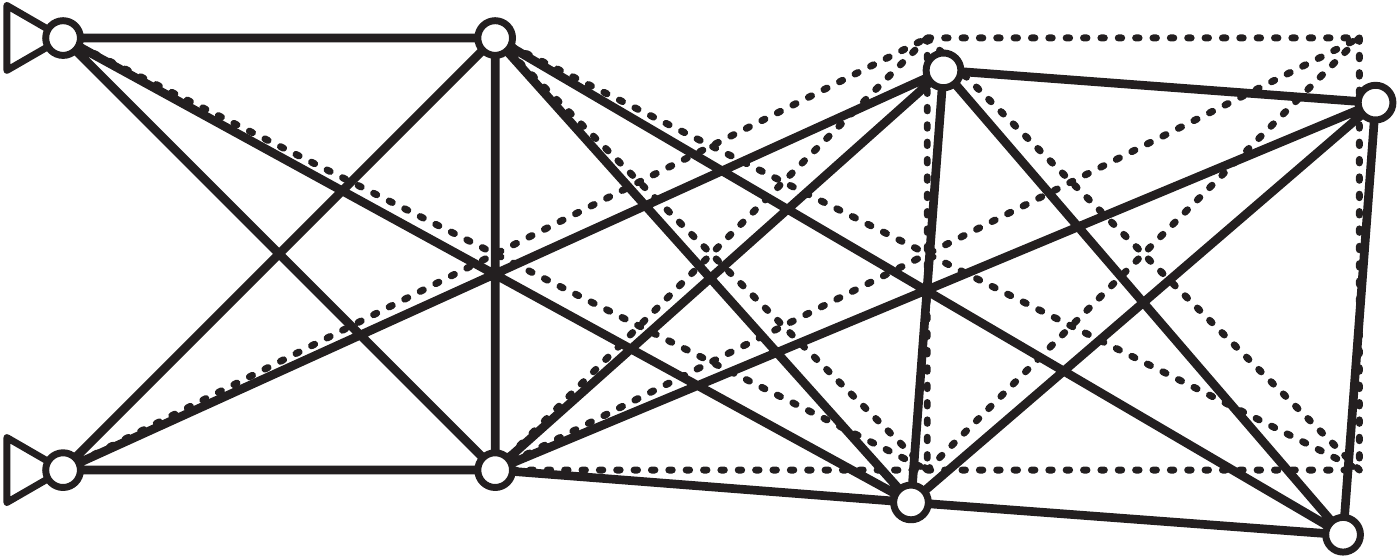}
    \caption{}
    \label{fig.enum2_alpha1}
  \end{subfigure}
  \par\bigskip
  \begin{subfigure}[b]{0.45\textwidth}
    \centering
    \includegraphics[scale=0.40]{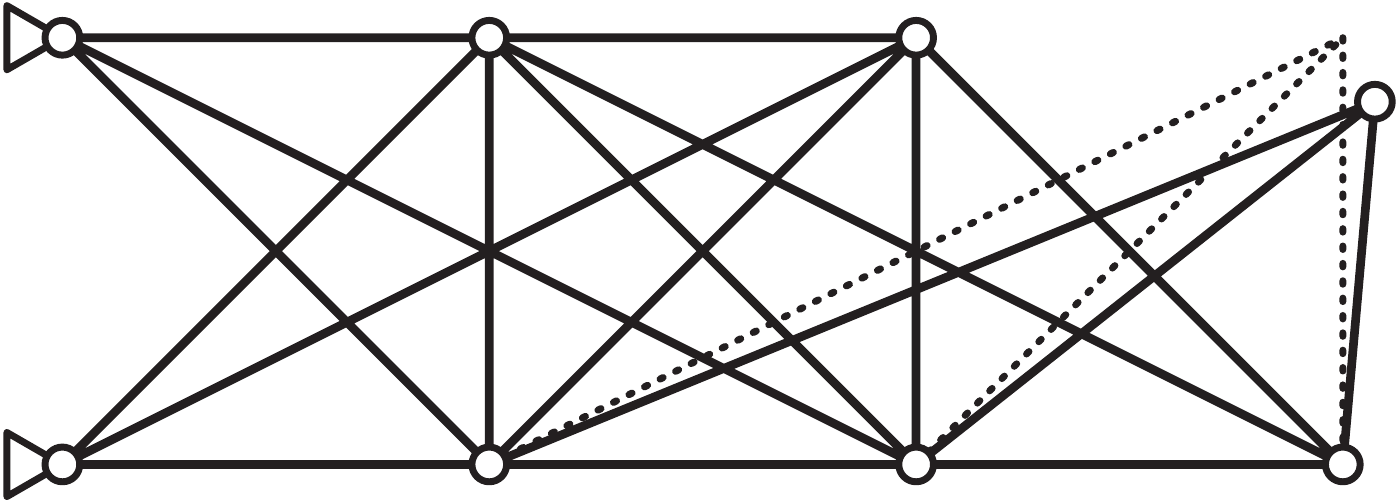}
    \caption{}
    \label{fig.enum3_alpha1}
  \end{subfigure}
  \hfill
  \begin{subfigure}[b]{0.45\textwidth}
    \centering
    \includegraphics[scale=0.40]{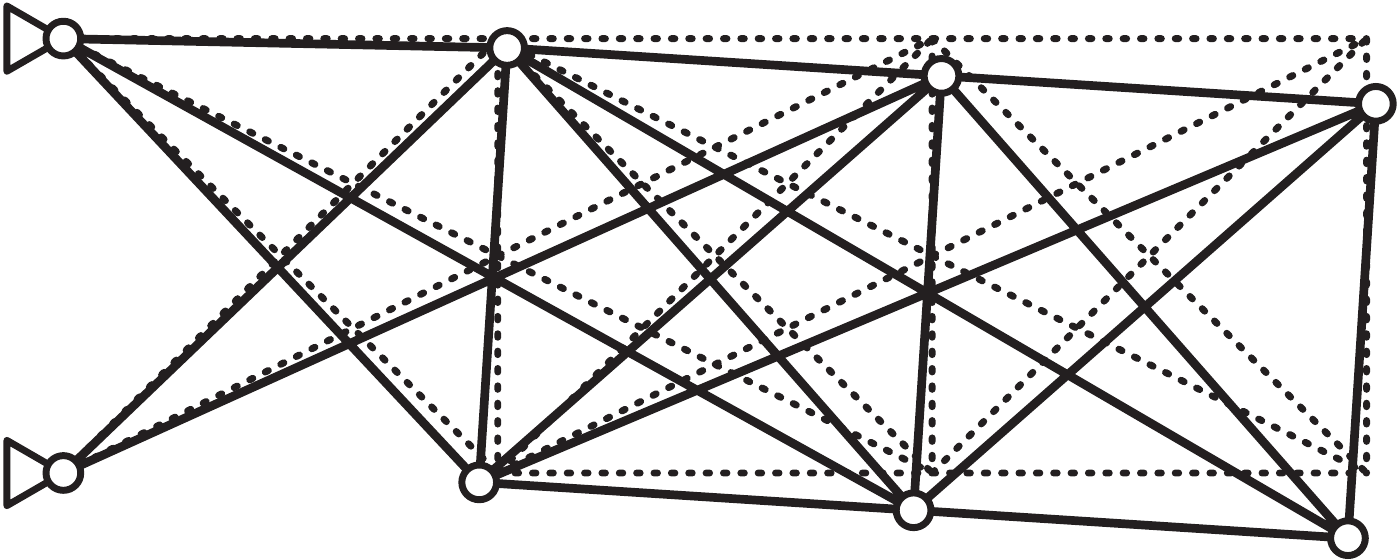}
    \caption{}
    \label{fig.enum4_alpha1}
  \end{subfigure}
  \par\bigskip
  \begin{subfigure}[b]{0.45\textwidth}
    \centering
    \includegraphics[scale=0.40]{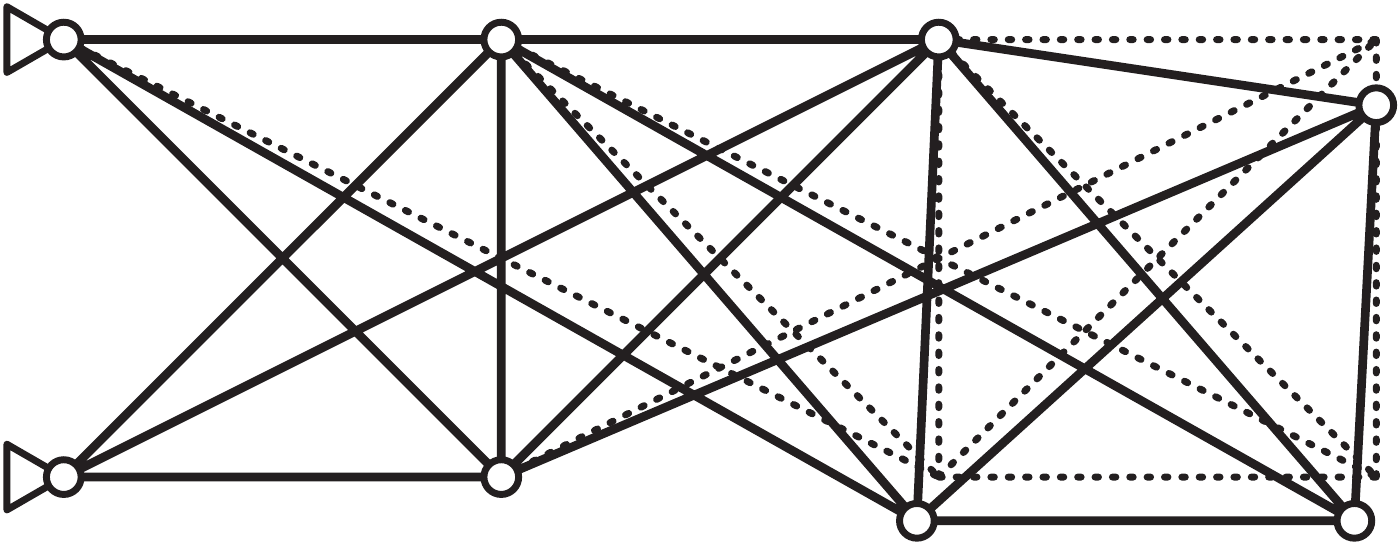}
    \caption{}
    \label{fig.enum5_alpha1}
  \end{subfigure}
  \hfill
  \begin{subfigure}[b]{0.45\textwidth}
    \centering
    \includegraphics[scale=0.40]{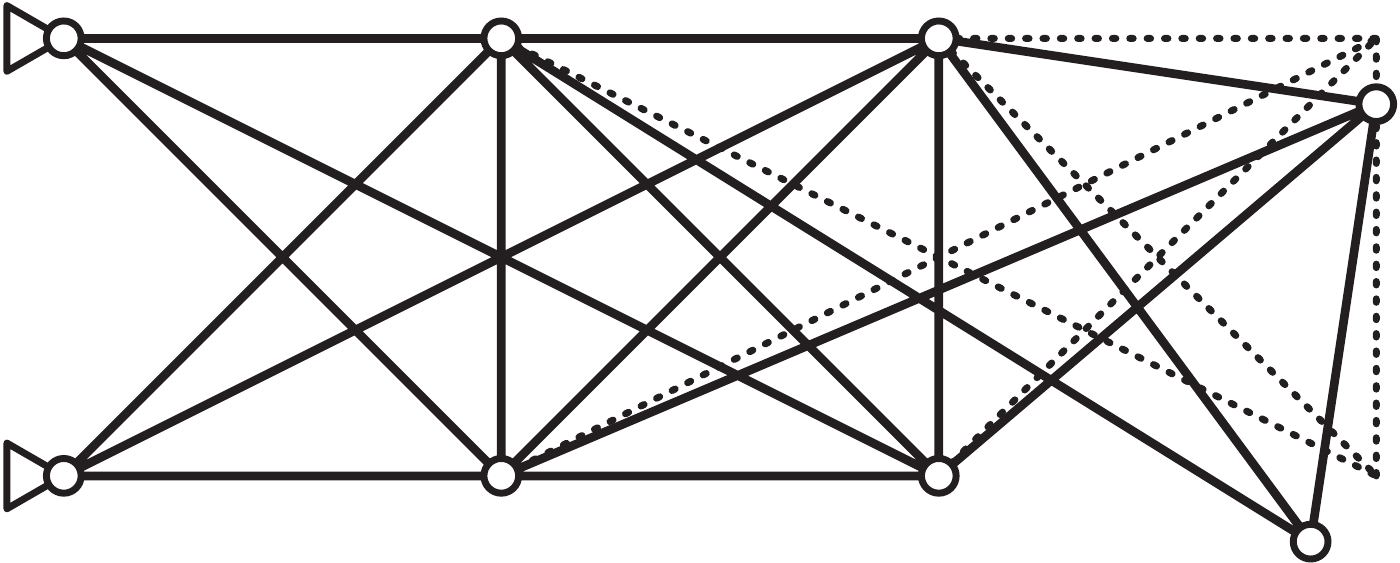}
    \caption{}
    \label{fig.enum6_alpha1}
  \end{subfigure}
  \par\bigskip
  \begin{subfigure}[b]{0.45\textwidth}
    \centering
    \includegraphics[scale=0.40]{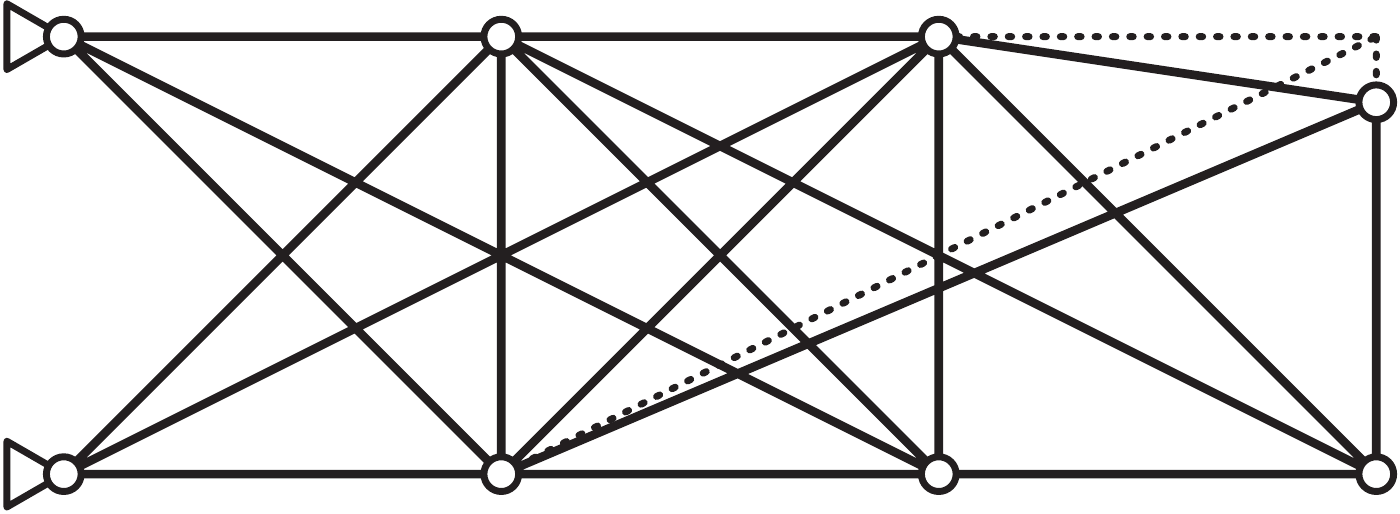}
    \caption{}
    \label{fig.enum15_alpha1}
  \end{subfigure}
  \hfill
  \begin{subfigure}[b]{0.45\textwidth}
  \end{subfigure}
  \caption{The worst-case scenarios for the solution with $\alpha=1$ in 
  example (I).}
  \label{fig.enum_alpha1}
\end{figure}

\refalg{alg.SQP} was implemented in MATLAB ver.\ 8.5.0 \citep{matlab}. 
In the algorithm, we need to solve some MILP problems at 
step~\ref{alg.step.filter.1} to evaluate the objective values at sample 
points, and solve a QP problem at step~\ref{alg.SQP.subproblem}. 
These MILP and QP problems are solved with 
CPLEX ver.~12.6.2 \citep{cplex}. 
Computation was carried out on 2.2{\,}GHz Intel Core i5-5200U 
processor with 8{\,}GB RAM. 

Consider the plane truss in \reffig{fig.m_19bar}, where $L=1\,\mathrm{m}$. 
The truss consists of $m=19$ members and has $d=12$ degrees of freedom 
of the displacements. 
As for the constant external load (i.e., $\bi{p}_{\rr{d}}$ in 
\eqref{P.lower.limit.analysis}), a horizontal force of 
$50\,\mathrm{kN}$ is applied at each of the rightmost nodes. 
A vertical force of $10\lambda\,\mathrm{kN}$ is applied at the upper 
rightmost node as the proportionally increasing load (i.e., 
$\lambda \bi{p}_{\rr{r}}$ in \eqref{P.lower.limit.analysis}). 
The yield stress is $\sigma_{\rr{y}}=200\,\mathrm{MPa}$. 

The initial point for \refalg{alg.SQP} is 
$\bi{x}_{0}=(1000,\dots,1000)^{\top}$ in $\mathrm{mm^{2}}$. 
The upper bound for the structural volume is given by 
$V=\bi{c}^{\top} \bi{x}_{0} = 2.6430 \times 10^{7}\,\mathrm{mm^{3}}$. 
The parameters for \refalg{alg.SQP} were chosen as 
$r=100\,\mathrm{mm^{2}}$, $r_{\min} = 10^{-4}\,\mathrm{mm^{2}}$, 
$\epsilon = 5 \times 10^{-4}\,\mathrm{mm^{2}}$, 
$\rho = 0.75$, 
$\eta = 0.01$, $\beta = 0.8$, $\tau_{\max}=50$, and 
$B_{0}$ is the identity matrix. 

%\begin{itemize}
%  \item 単位の換算
%%cs = 10 * ones(nm,1);
%%sigmaYield = 20.0;
%%qYield = sigmaYield * cs;
%%vec_f.D(10) = -50.0;
%%vec_f.D(12) = -50.0;
%%vec_f.R(9)  =  10;
%  \begin{itemize}
%    \item 実装：
%          $cs_{i} = 10$, $\sigma_{\rr{y}}=20$, 
%          $q_{\rr{y}i} = cs_{i} \times \sigma_{\rr{y}} = 20 \times 10 = 200$, 
%          $f_{\rr{d}j}=50$, $f_{\rr{r}j}=10$ 
%    \item 実際の単位：
%          $x_{i} = 10\,\mathrm{cm^{2}} = 1000\,\mathrm{mm^{2}}$, 
%          $\sigma_{\rr{y}} = 200\,\mathrm{MPa} = 200\,\mathrm{N/mm^{2}}$, 
%          $q_{\rr{y}i} 
%          = (1000\,\mathrm{mm^{2}}) \times (200\,\mathrm{N/mm^{2}})
%          = 200\,\mathrm{kN}$, 
%          $f_{\rr{d}j}=50\,\mathrm{kN}$, $f_{\rr{r}j}=10\,\mathrm{kN}$ 
%  \end{itemize}
%
%  \item $\bi{x}_{0} = (10,\dots,10)^{\top}$. 
%        体積の上限値は，$\bi{x}_{0}$ の体積とする．
%  \begin{itemize}
%    \item つまり，$264.2955$（単位は $\mathrm{cm^{2}} \times \mathrm{m}$）．
%    \item $10 \times 10 \times 1000 = 10^{5}$ 倍すれば，
%          $\mathrm{mm^{3}}$ 単位になる $\to$ 
%          $264.2955 \times 10^{5}\,\mathrm{mm^{3}}
%          = 2.642955 \times 10^{7}\,\mathrm{mm^{3}}$
%  \end{itemize}
%  \item $r$ の初期値は $0.1 \bi{x}_{0} = 1$，
%        $r_{\min} = 10^{-6}$. 
%        $\epsilon = 5 \times 10^{-6}$. 
%  \item $\eta = 0.01$, $\beta = 0.8$, 
%        $\tau_{\max}$ は $\alpha$ の最小値が $10^{-5}$ になるように選
%        ぶ（つまり，$\tau_{\max}=51$）．
%  \item $\rho = 0.75$
%  \item $B_{0}$ は単位行列．
%  \item $f_{\rr{d}}=-50.0$ （右端の二つの頂点に左向き）．
%        $f_{\rr{r}}=10$ （右端上の頂点に下向き）．
%        降伏軸力は，$(\text{断面積}) \times 20$ として計算．
%\end{itemize}

For $\alpha=1$, the solution obtained by \refalg{alg.SQP} is 
shown in \reffig{fig.optimal_alpha_1}, where 
the width of each member is proportional to its cross-sectional area. 
The algorithm terminates after solving $376$ QP problems. 
The number of MILP problems solved for the objective function 
evaluations is $3699$. 
The obtained solution satisfies the termination 
condition $\| \bi{d}_{k} \| < \epsilon$ with a small value of the stencil 
radius, $r=4.2 \times 10^{-3}\,\mathrm{mm^{2}}$. 
%Hence, we may expect that the simplex gradient at the convergent solution 
%is close to the true gradient (if it exists). 
%This might guarantee the local optimality of the obtained solution; see 
%\refrem{rem.optimality}. 
The worst-case limit load factor of the obtained solution is $14.4979$, 
while that of the initial design is $6.7187$. 
The worst-case scenarios for the obtained solution, as well as the corresponding 
collapse modes, are collected in \reffig{fig.enum_alpha1}, where the 
damaged members are removed from the figures. 
It is emphasized that the limit load factors of these seven scenarios 
are all equal to the objective value, $14.4979$. 
In contrast, the worst-case scenario for the initial design is only the 
one shown in \reffig{fig.enum4_alpha1}. 
It seems to be natural that a local optimal solution of a redundancy 
optimization problem has multiple worst-case scenarios in general. 
Namely, multiplicity of worst-case scenarios means that, 
if we attempt to increase the limit load factor corresponding to 
a certain scenario, then the one corresponding to another scenario 
decreases. 

\begin{figure}[tp]
  %%% C:\doc\derivative_free\redundancy\eva5\opt_design.m
  \centering
  \begin{subfigure}[b]{0.45\textwidth}
    \centering
    \includegraphics[scale=0.40]{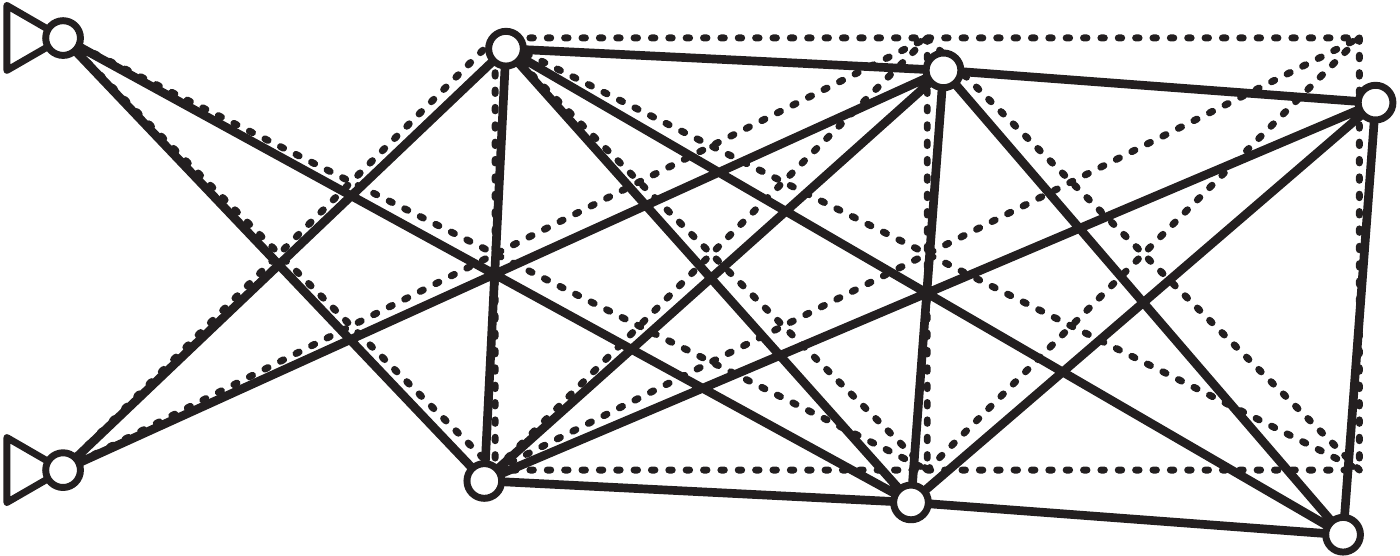}
    \caption{}
    \label{fig.enum1_4_alpha2}
  \end{subfigure}
  \hfill
  \begin{subfigure}[b]{0.45\textwidth}
    \centering
    \includegraphics[scale=0.40]{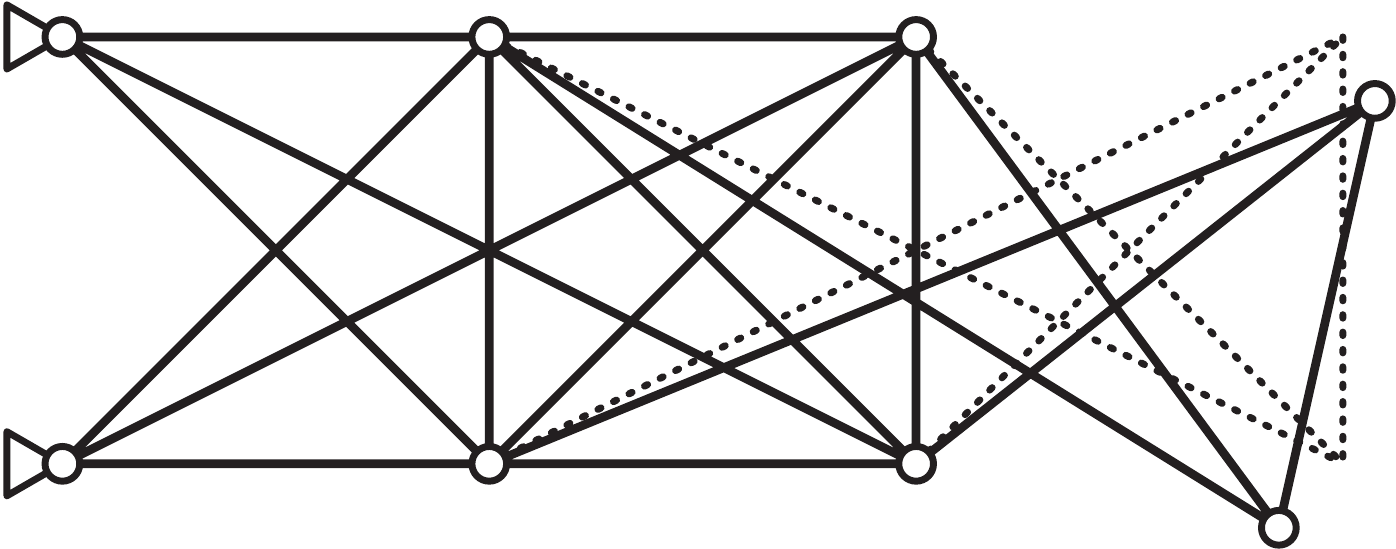}
    \caption{}
    \label{fig.enum3_6_alpha2}
  \end{subfigure}
  \par\bigskip
  \begin{subfigure}[b]{0.45\textwidth}
    \centering
    \includegraphics[scale=0.40]{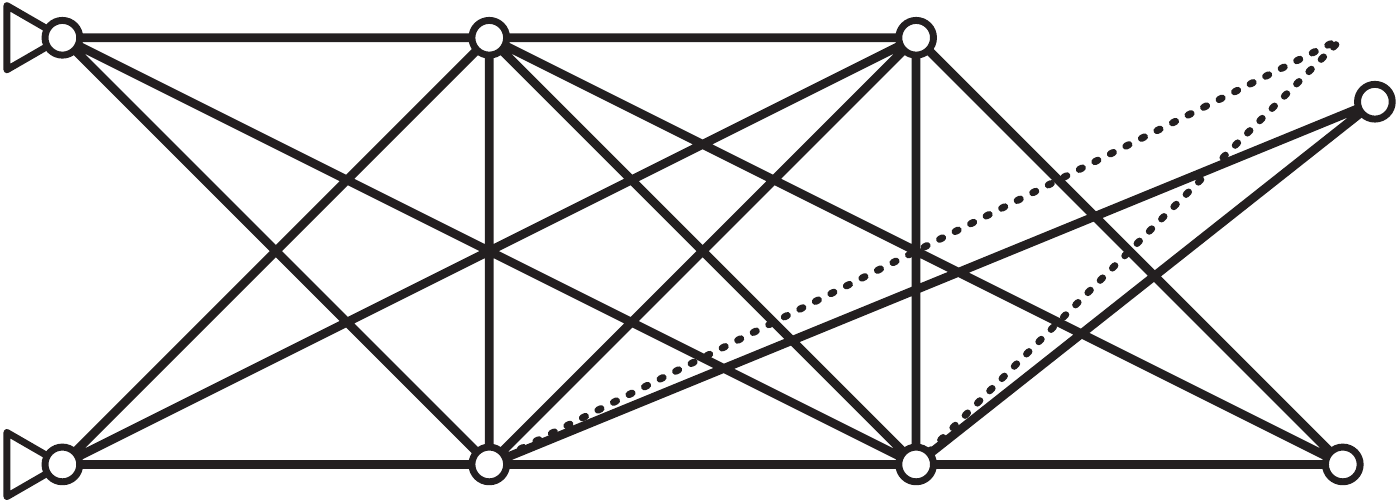}
    \caption{}
    \label{fig.enum3_9_alpha2}
  \end{subfigure}
  \hfill
  \begin{subfigure}[b]{0.45\textwidth}
    \centering
    \includegraphics[scale=0.40]{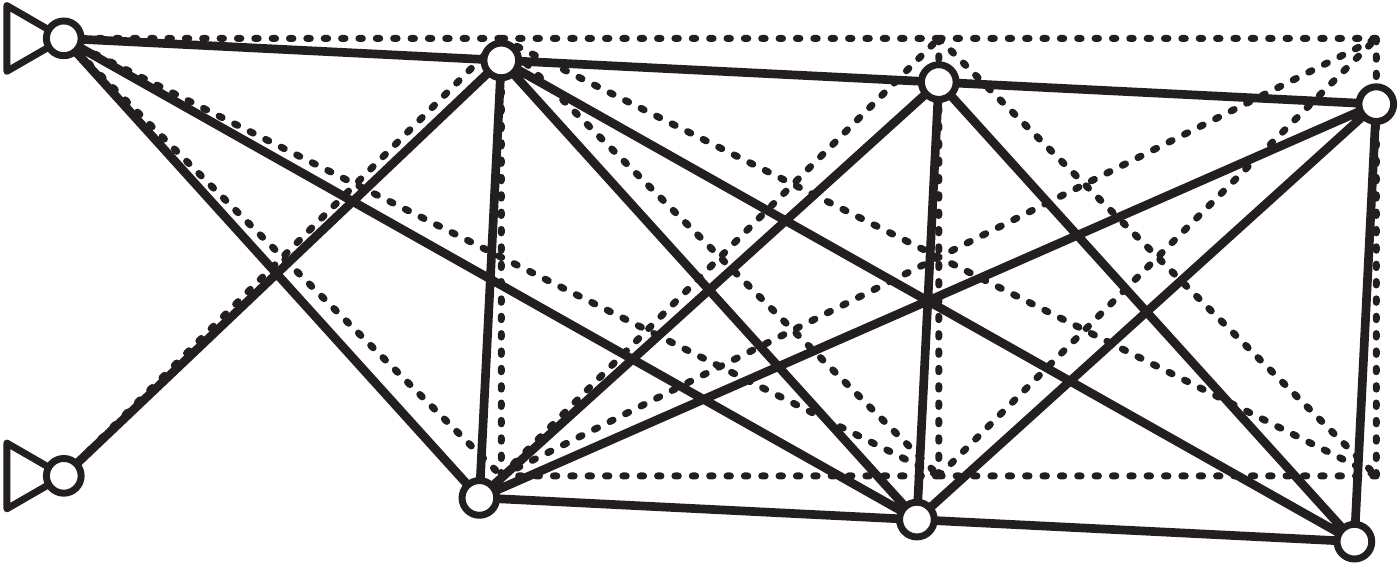}
    \caption{}
    \label{fig.enum4_18_alpha2}
  \end{subfigure}
  \par\bigskip
  \begin{subfigure}[b]{0.45\textwidth}
    \centering
    \includegraphics[scale=0.40]{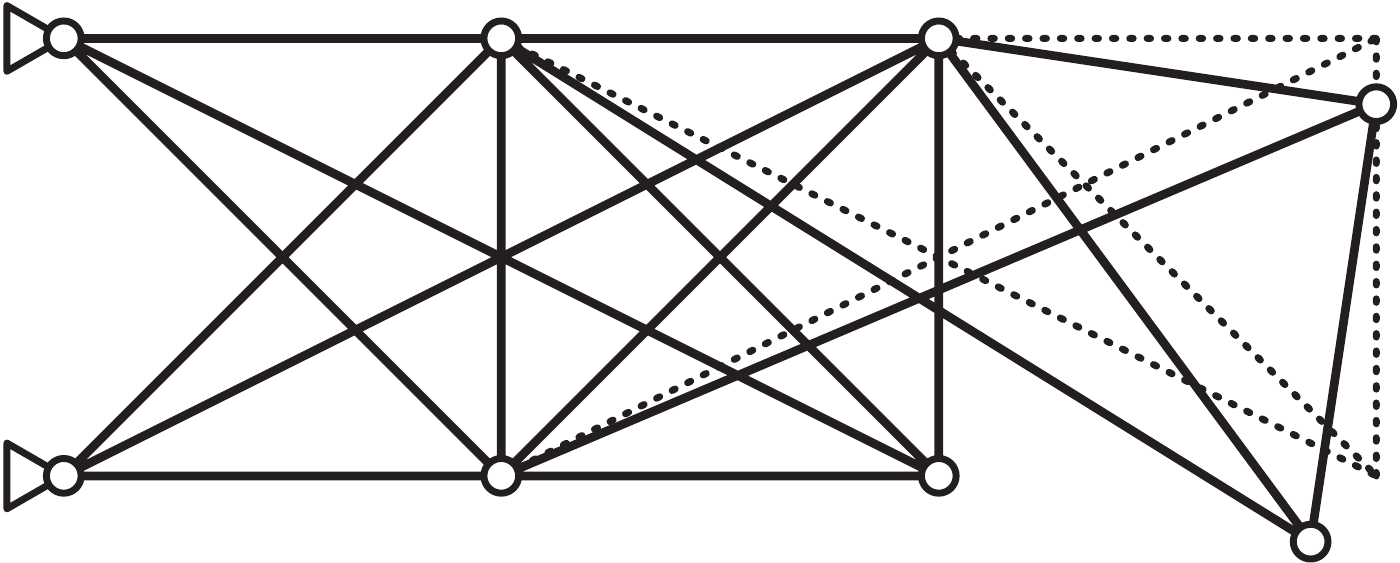}
    \caption{}
    \label{fig.enum6_15_alpha2}
  \end{subfigure}
  \hfill
  \begin{subfigure}[b]{0.45\textwidth}
    \centering
    \includegraphics[scale=0.40]{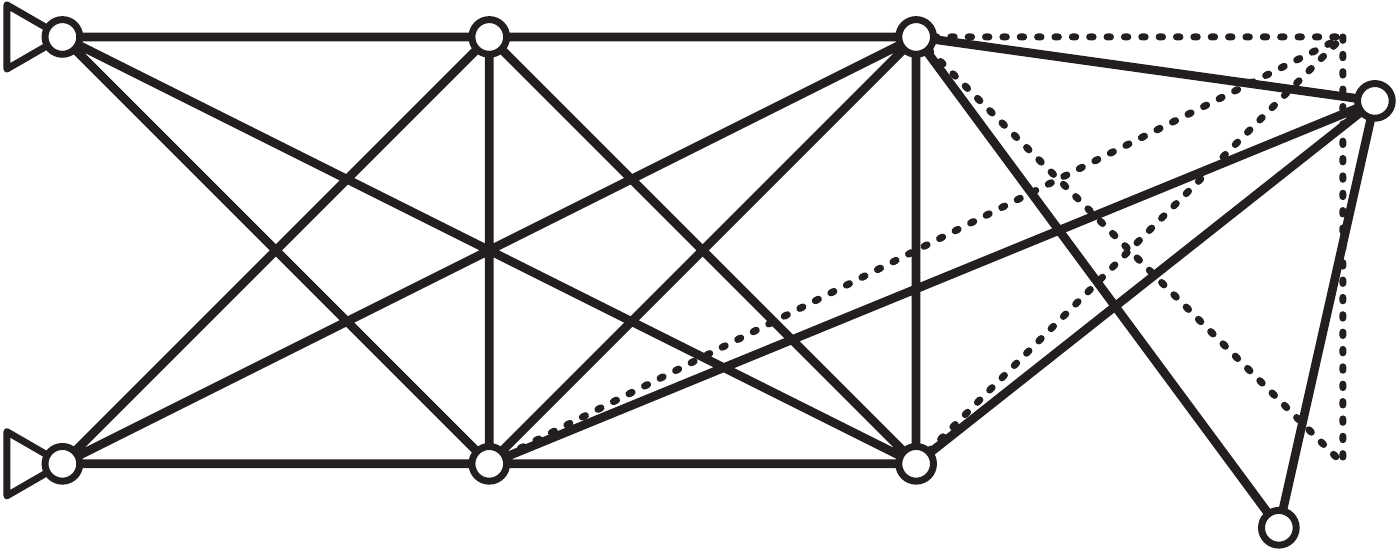}
    \caption{}
    \label{fig.enum6_17_alpha2}
  \end{subfigure}
  \par\bigskip
  \begin{subfigure}[b]{0.45\textwidth}
    \centering
    \includegraphics[scale=0.40]{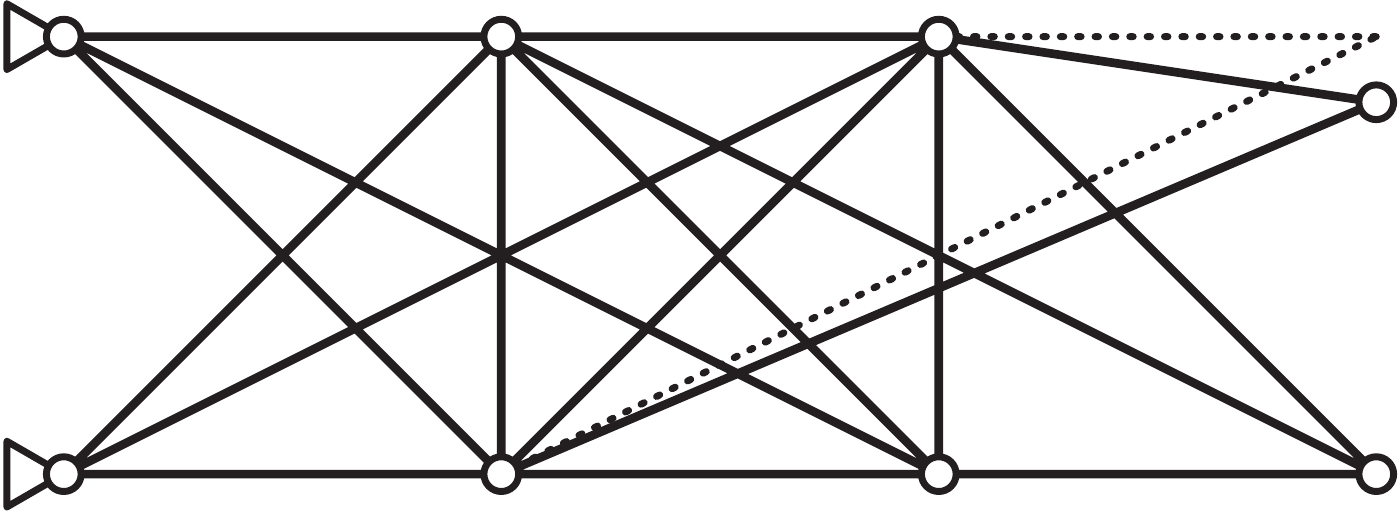}
    \caption{}
    \label{fig.enum9_15_alpha2}
  \end{subfigure}
  \hfill
  \begin{subfigure}[b]{0.45\textwidth}
    \centering
    \includegraphics[scale=0.40]{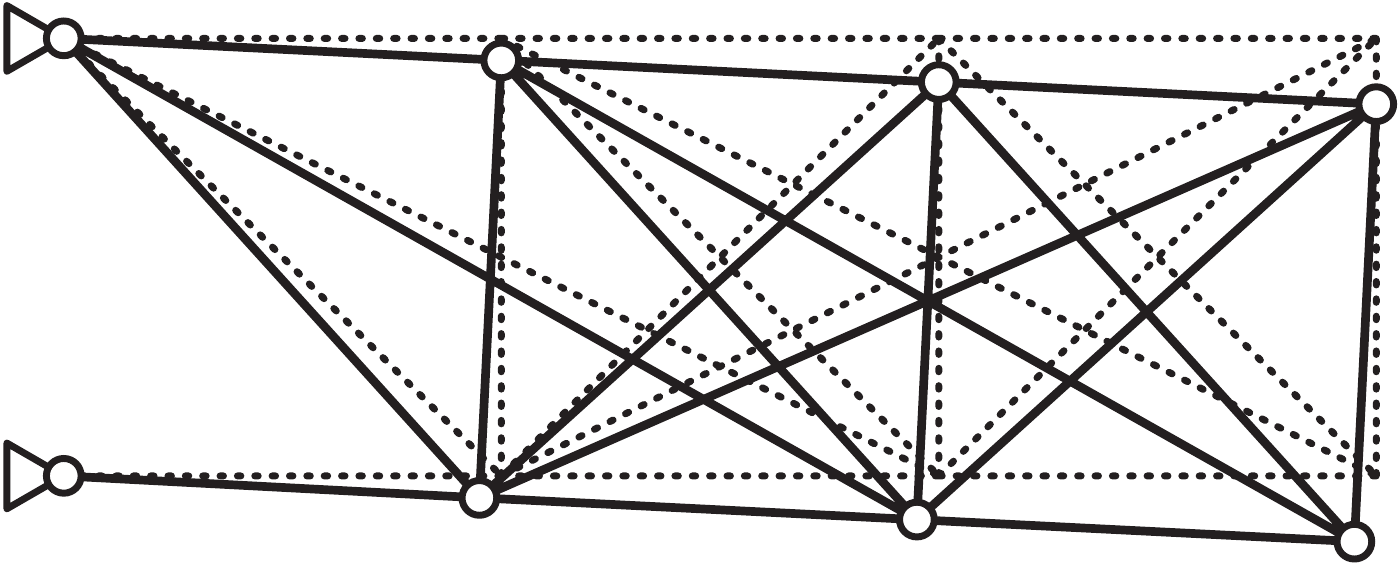}
    \caption{}
    \label{fig.enum13_18_alpha2}
  \end{subfigure}
  \par\bigskip
  \begin{subfigure}[b]{0.45\textwidth}
    \centering
    \includegraphics[scale=0.40]{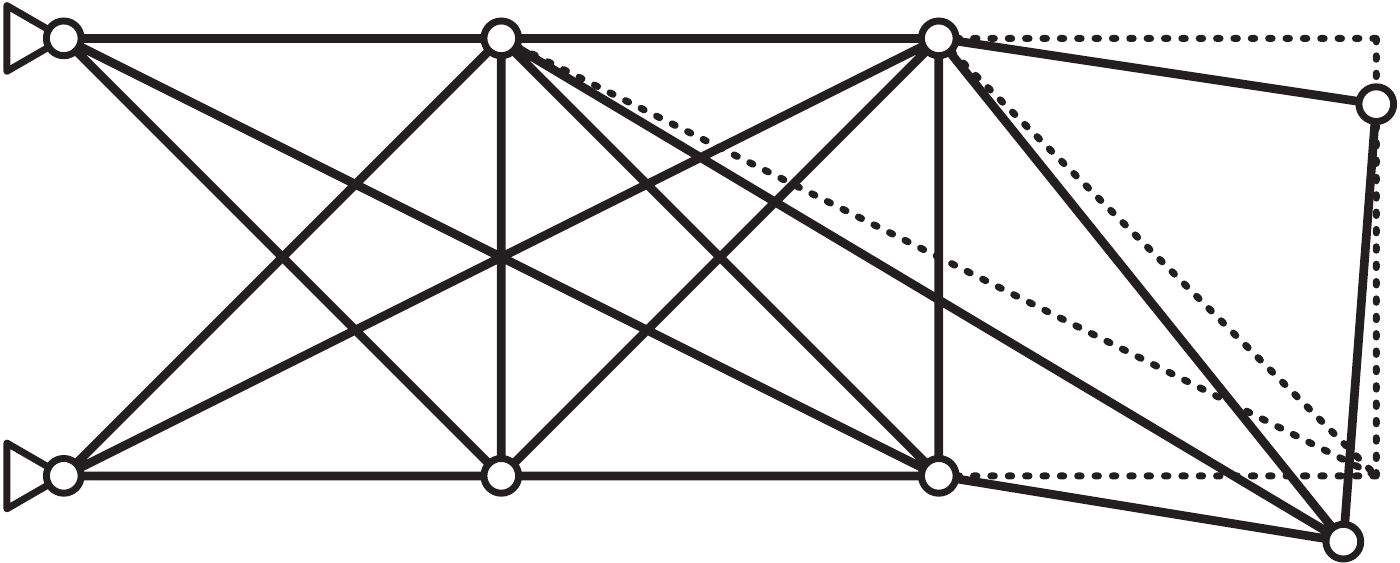}
    \caption{}
    \label{fig.enum15_19_alpha2}
  \end{subfigure}
  \hfill
  \begin{subfigure}[b]{0.45\textwidth}
  \end{subfigure}
  \caption{The worst-case scenarios for the solution with $\alpha=2$ 
  in example (I).}
  \label{fig.enum_alpha2}
\end{figure}

We next examine $\alpha=2$. 
The solution obtained by \refalg{alg.SQP} is shown in 
\reffig{fig.optimal_alpha_2}. 
The algorithm terminates after solving $327$ QP problems and $3326$ MILP 
problems. 
The obtained solution satisfies the termination 
condition $\| \bi{d}_{k} \| < \epsilon$ with a small value of the stencil 
radius, $r=2.4\times 10^{-3}\,\mathrm{mm^{2}}$. 
%Therefore, the obtained solution might be considered to satisfy the 
%local optimality condition. 
The worst-case limit load factor of the obtained solution is 
$6.5509$, while that of the initial design is $3.0474$. 
\reffig{fig.enum_alpha2} collects the worst-case scenarios for the obtained 
solution. 
Namely, the multiplicity of the worst-case scenarios is $9$. 
In contrast, the worst-case scenario for the initial design is unique and 
is the one shown in \reffig{fig.enum4_18_alpha2}. 
It is worth noting that the objective function is not differentiable in 
general at a point having multiple worst-case scenarios. 
Nevertheless, the proposed algorithm could find a solution with large 
multiplicity. 
%\OMIT{このように implicit filtering が nonsmooth optimization に有用なこと
%は，\cite{CLLLR15,CDV08,CV12,RS13} などに述べられている?}

The ground structure in \reffig{fig.m_19bar} becomes unstable if a 
particular set of three members is removed. 
This means that the redundancy optimization problem, 
\eqref{P.redundancy.optimization}, loses its meaning for 
$\alpha \ge 3$. 
This is because we assume that a damaged structural component is 
completely absent from the structure. 
In contrast, if we adopt a nonzero degree of damage, $\gamma$, as 
discussed in \refrem{rem.partial.deficiency}, then problem 
\eqref{P.redundancy.optimization} has meaning even for $\alpha \ge 3$. 

When we set $\alpha=0$, structural degradation is not 
considered, and the redundancy optimization of the limit load factor 
reverts to the conventional limit design (the optimal plastic design). 
The optimal solution of the limit design problem, obtained by linear 
programming, is shown in \reffig{fig.optimal_alpha_0}. 
%The limit load factor of this solution is $34.9061$. 
It is worth noting that this is a statically determinate truss. 
Hence, the truss becomes unstable (kinematically indeterminate) if 
any single member is removed. 
Therefore, the truss has no redundancy; more precisely, the strong 
redundancy defined by \cite{KBh11} is equal to zero. 

\begin{figure}[tp]
  \centering
  \begin{subfigure}[b]{0.55\textwidth}
    \centering
    \scalebox{0.80}{
    \includegraphics{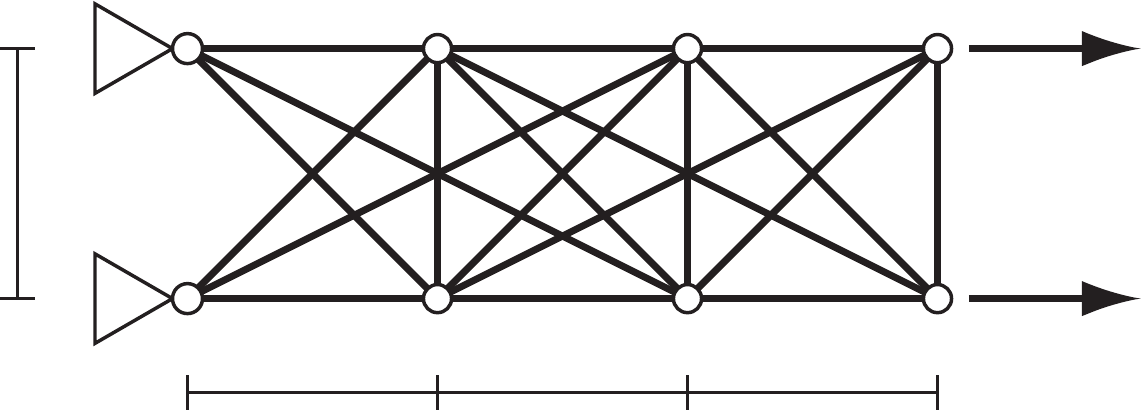}
    \begin{picture}(0,0)
      \put(-288,20){
      \put(-52,45){{\large $L$}}
      \put(245,94){{\large $\lambda p_{\rr{r}}$}}
      \put(245,-2){{\large $\lambda p_{\rr{r}}$}}
      \put(42,-28){{\large $L$}}
      \put(114,-28){{\large $L$}}
      \put(186,-28){{\large $L$}}
      }
    \end{picture}
    }
    \caption{}
    \label{fig.m_19bar_2}
  \end{subfigure}
  \par\bigskip\bigskip
  \begin{subfigure}[b]{0.55\textwidth}
    \centering
    %%% C:\doc\derivative_free\redundancy\eva7_2\opt_design.m
    \includegraphics[scale=0.50]{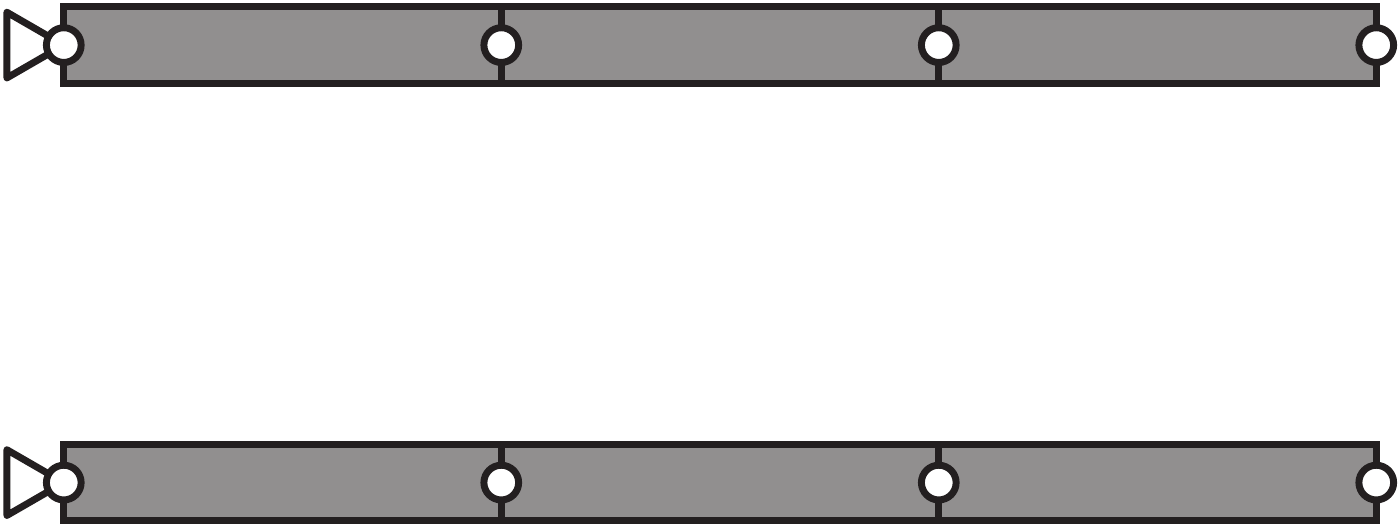}
    \caption{}
    \label{fig.ex2_optimal_alpha_0}
  \end{subfigure}
  \caption{Example (II). 
  \subref{fig.m_19bar_2}~Problem setting; and 
  \subref{fig.ex2_optimal_alpha_0}~the optimal solution without 
  considering redundancy (i.e., $\alpha=0$). 
  }
\end{figure}

\begin{figure}[tp]
  \centering
  \begin{subfigure}[b]{0.55\textwidth}
    \centering
    %%% C:\doc\derivative_free\redundancy\eva9\opt_design.m
    \includegraphics[scale=0.50]{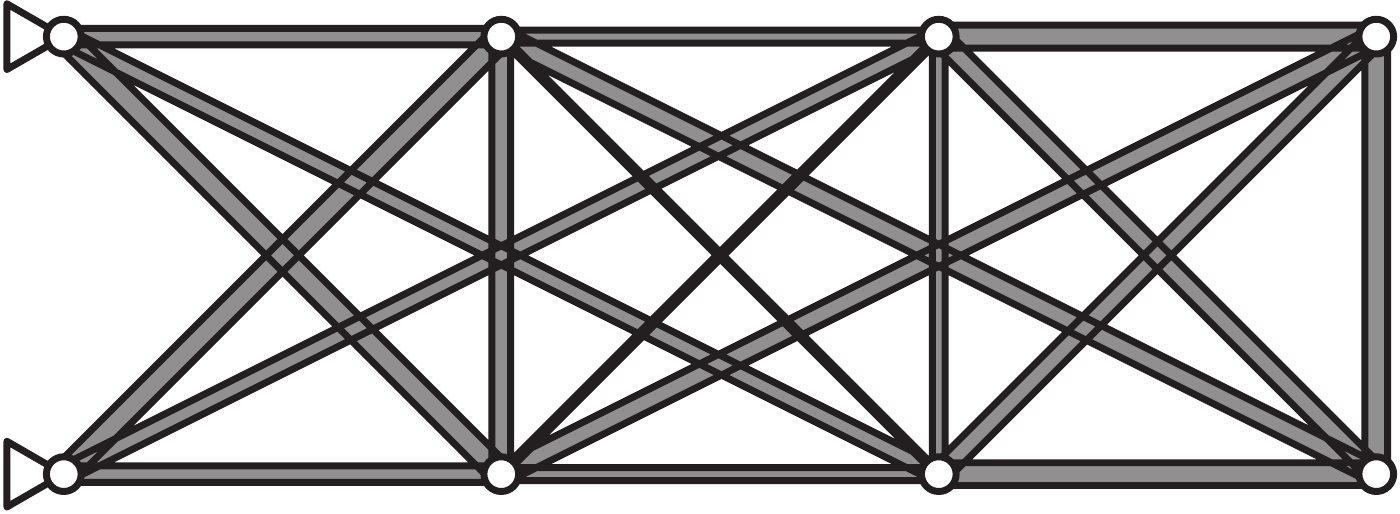}
    \caption{}
    \label{fig.ex2_optimal_alpha_1}
  \end{subfigure}
  \par\bigskip\bigskip
  \begin{subfigure}[b]{0.55\textwidth}
    \centering
    %%% C:\doc\derivative_free\redundancy\eva8\opt_design.m
    \includegraphics[scale=0.50]{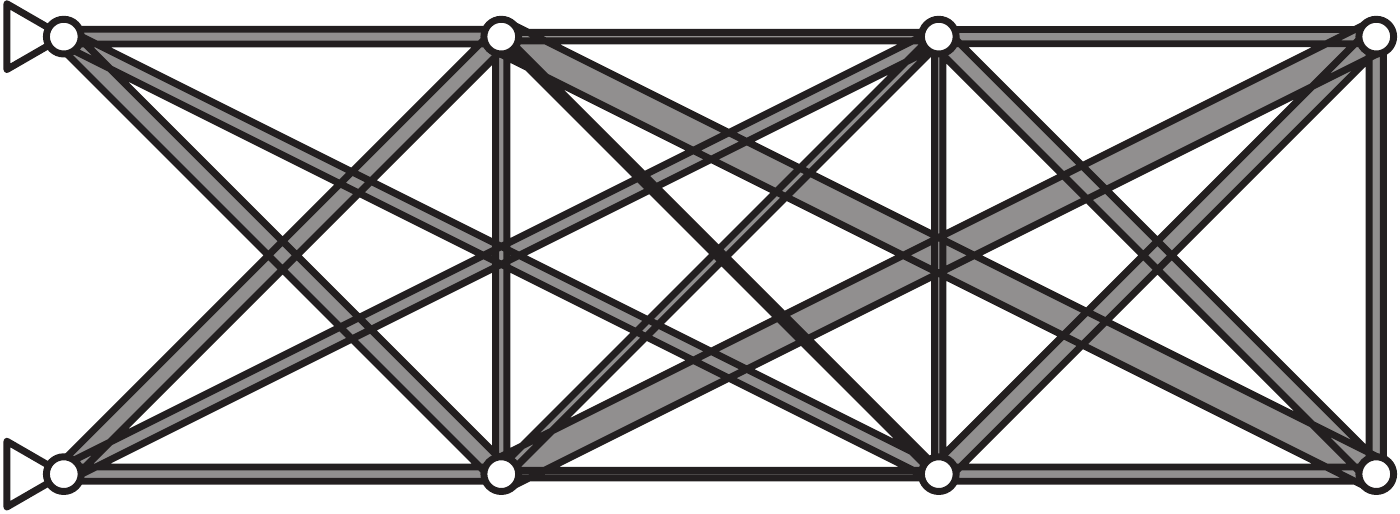}
    \caption{}
    \label{fig.ex2_optimal_alpha_2}
  \end{subfigure}
  \caption{The solutions obtained in example (II).
  \subref{fig.ex2_optimal_alpha_1}~$\alpha=1$; and 
  \subref{fig.ex2_optimal_alpha_2}~$\alpha=2$. }
\end{figure}

\begin{figure}[tp]
  %%% C:\doc\derivative_free\redundancy\eva9\opt_design.m
  \centering
  \begin{subfigure}[b]{0.45\textwidth}
    \centering
    \includegraphics[scale=0.40]{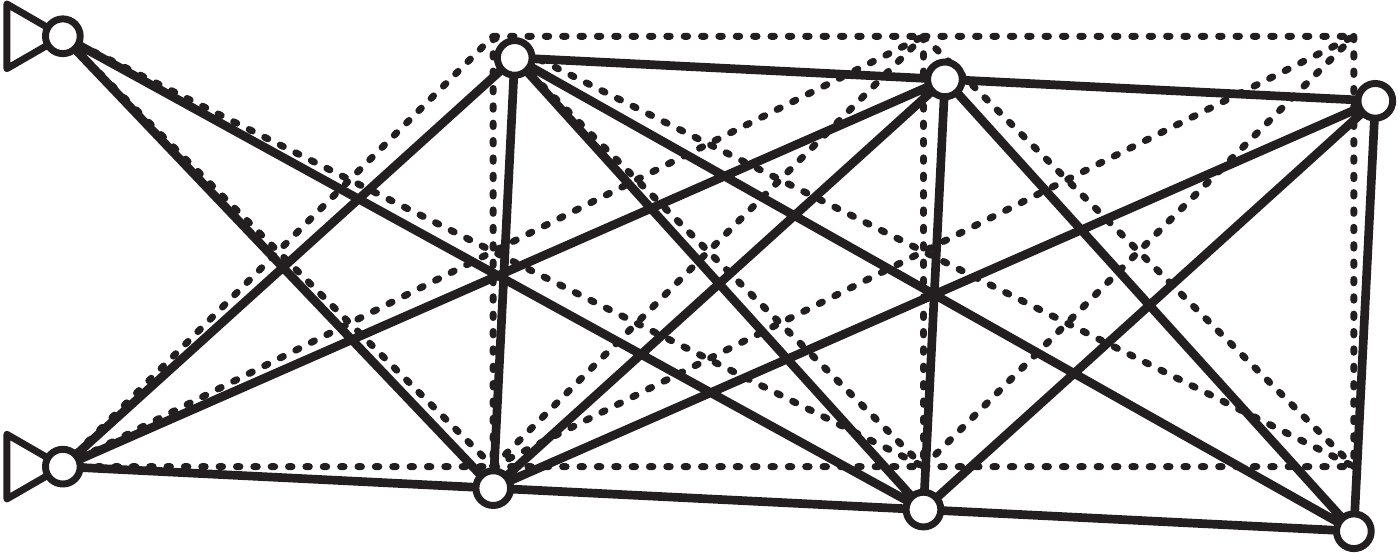}
    %\par
    %\includegraphics[scale=0.40]{ex2_enum4_alpha1.pdf}
    \caption{}
    \label{fig.ex2_enum1_alpha1}
  \end{subfigure}
  \hfill
  \begin{subfigure}[b]{0.45\textwidth}
    \centering
    \includegraphics[scale=0.40]{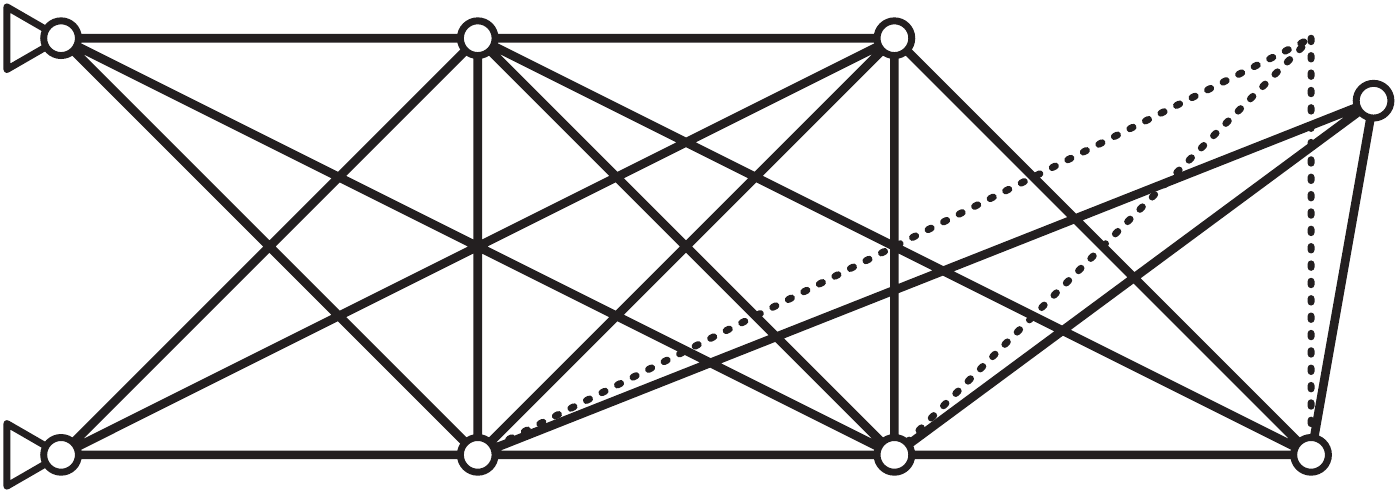}
    %\par
    %\includegraphics[scale=0.40]{ex2_enum6_alpha1.pdf}
    \caption{}
    \label{fig.ex2_enum3_alpha1}
  \end{subfigure}
  \par\bigskip
  \begin{subfigure}[b]{0.45\textwidth}
    \centering
    \includegraphics[scale=0.40]{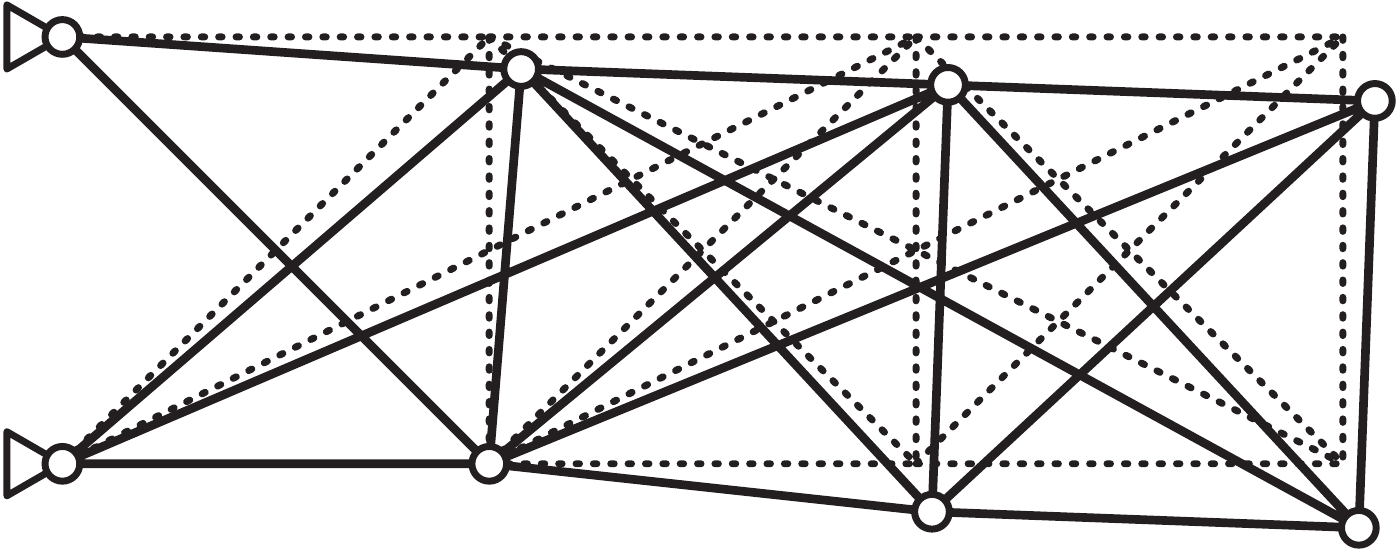}
    %\par
    %\includegraphics[scale=0.40]{ex2_enum18_alpha1.pdf}
    \caption{}
    \label{fig.ex2_enum16_alpha1}
  \end{subfigure}
  \hfill
  \begin{subfigure}[b]{0.45\textwidth}
    \centering
    \includegraphics[scale=0.40]{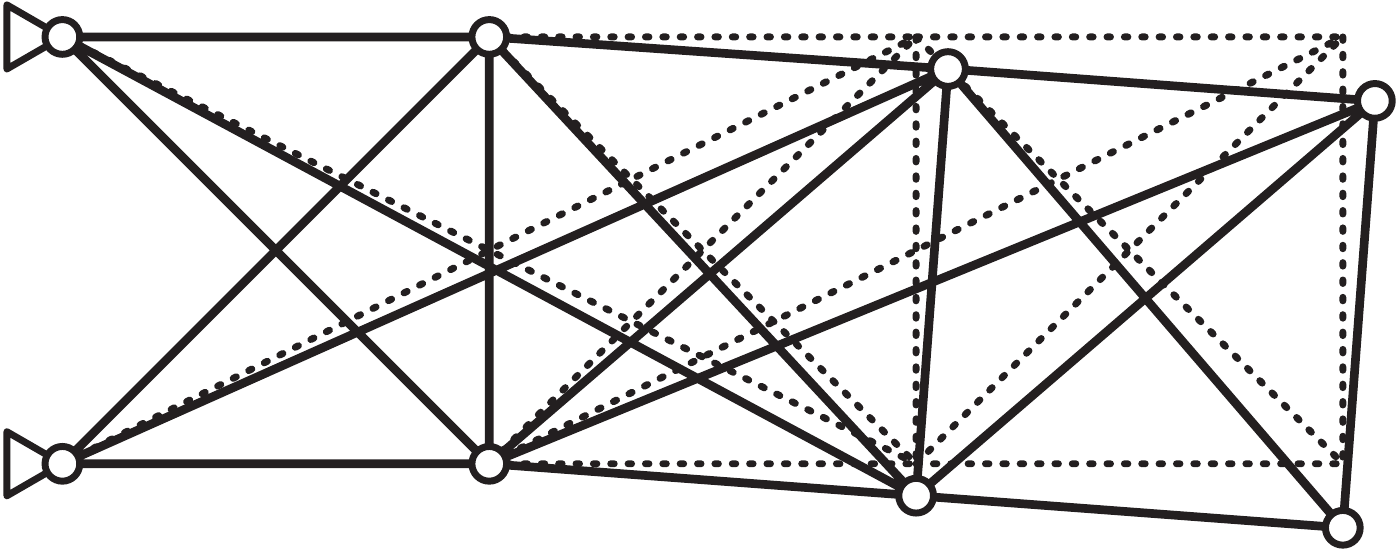}
    %\par
    %\includegraphics[scale=0.40]{ex2_enum19_alpha1.pdf}
    \caption{}
    \label{fig.ex2_enum17_alpha1}
  \end{subfigure}
  \par\bigskip
  \begin{subfigure}[b]{0.45\textwidth}
    \centering
    \includegraphics[scale=0.40]{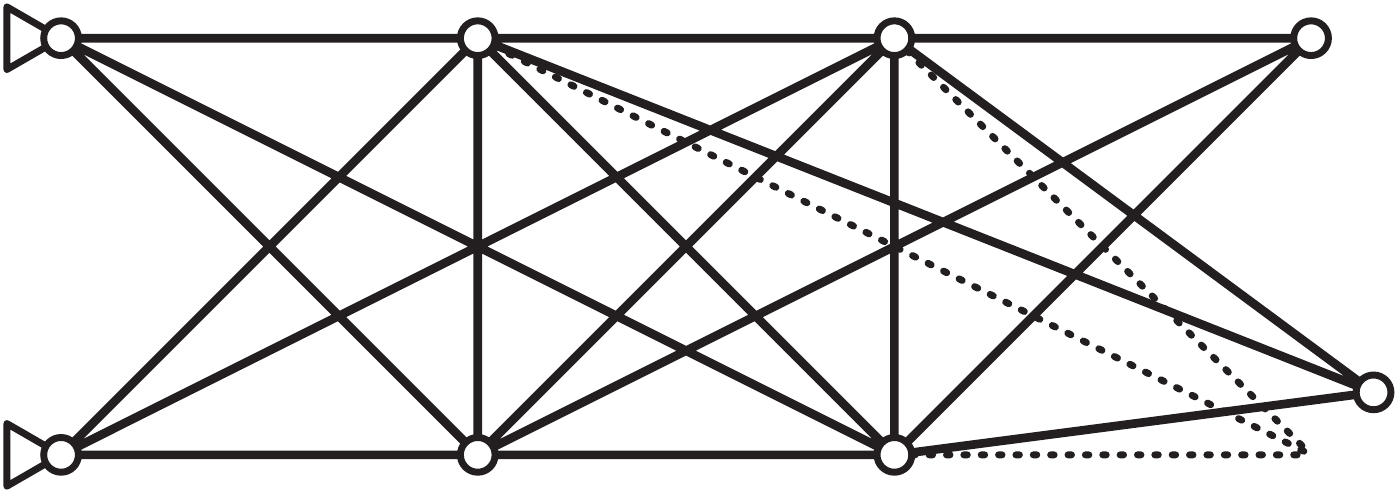}
    \caption{}
    \label{fig.ex2_enum9_alpha1}
  \end{subfigure}
  \hfill
  \begin{subfigure}[b]{0.45\textwidth}
  \end{subfigure}
  \caption{Some of the worst-case scenarios for the solution with 
  $\alpha=1$ in example (II). 
  The damage scenarios obtained by reflecting the scenarios in 
  Figures~\ref{fig.ex2_enum1_alpha1}, 
  \subref{fig.ex2_enum3_alpha1}, \subref{fig.ex2_enum16_alpha1}, and 
  \subref{fig.ex2_enum17_alpha1} across the axis of symmetry are also 
  the worst-case scenarios. 
  }
  \label{fig.ex2_enum_alpha1}
\end{figure}

\begin{figure}[tp]
  %%% C:\doc\derivative_free\redundancy\eva8\opt_design.m
  \centering
  \begin{subfigure}[b]{0.45\textwidth}
    \centering
    \includegraphics[scale=0.40]{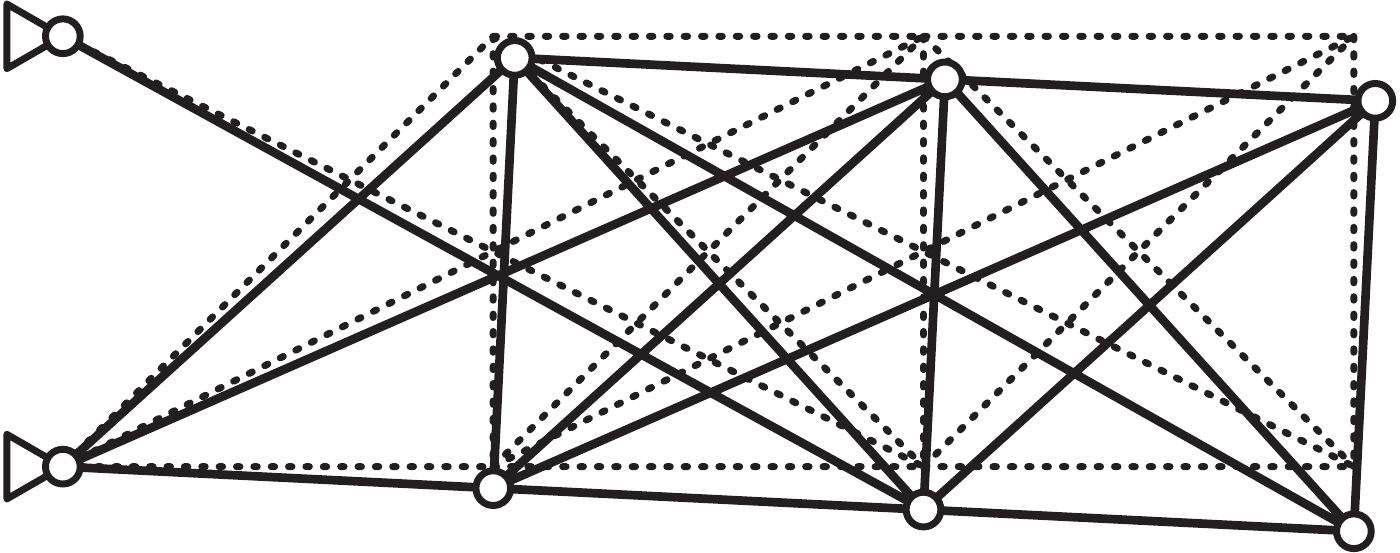}
    %\par
    %\includegraphics[scale=0.40]{ex2_enum4_13_alpha2.pdf}
    \caption{}
    \label{fig.ex2_enum1_10_alpha2}
  \end{subfigure}
  \hfill
  \begin{subfigure}[b]{0.45\textwidth}
    \centering
    \includegraphics[scale=0.40]{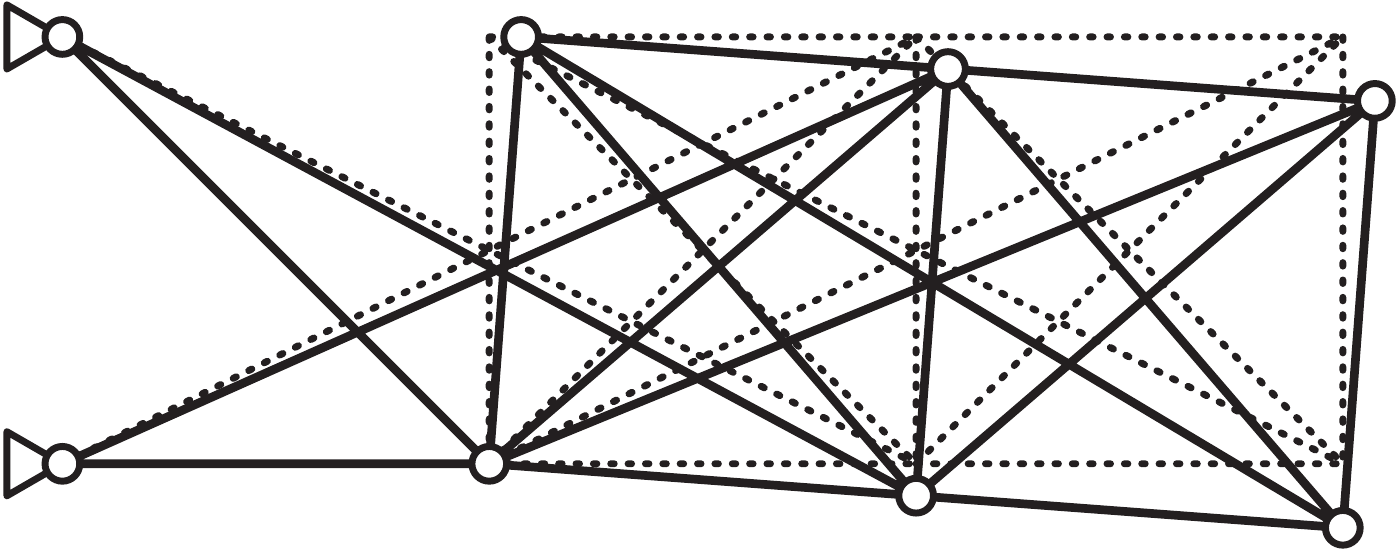}
    %\par
    %\includegraphics[scale=0.40]{ex2_enum4_10_alpha2.pdf}
    \caption{}
    \label{fig.ex2_enum1_13_alpha2}
  \end{subfigure}
  \par\bigskip
  \begin{subfigure}[b]{0.45\textwidth}
    \centering
    \includegraphics[scale=0.40]{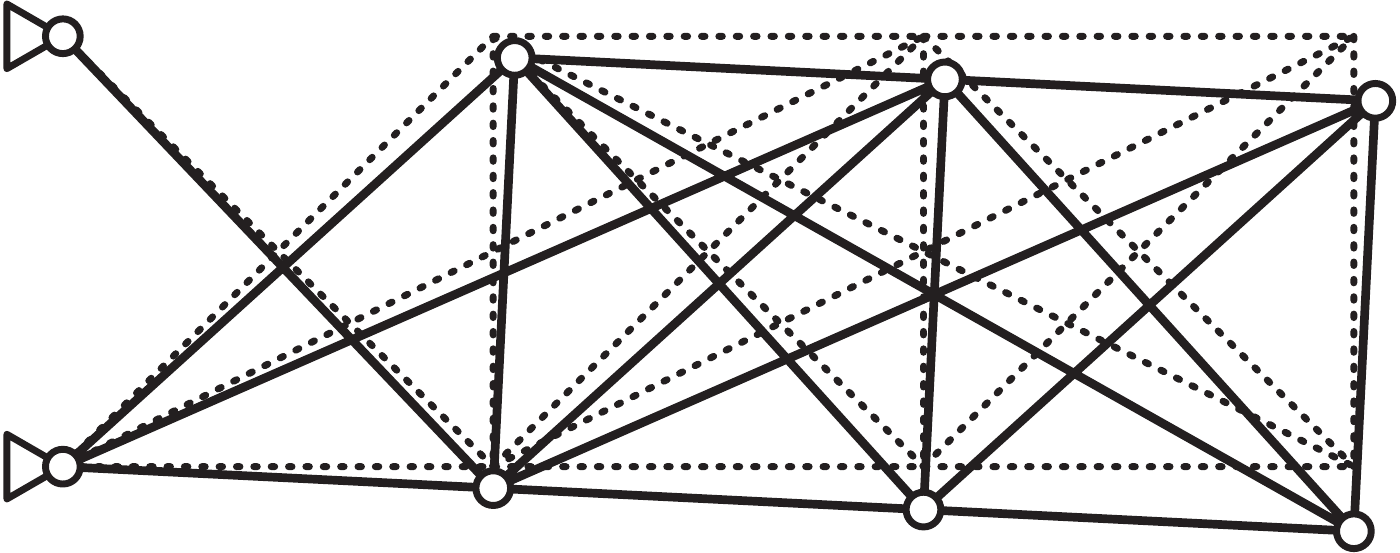}
    %\par
    %\includegraphics[scale=0.40]{ex2_enum4_18_alpha2.pdf}
    \caption{}
    \label{fig.ex2_enum1_16_alpha2}
  \end{subfigure}
  \hfill
  \begin{subfigure}[b]{0.45\textwidth}
    \centering
    \includegraphics[scale=0.40]{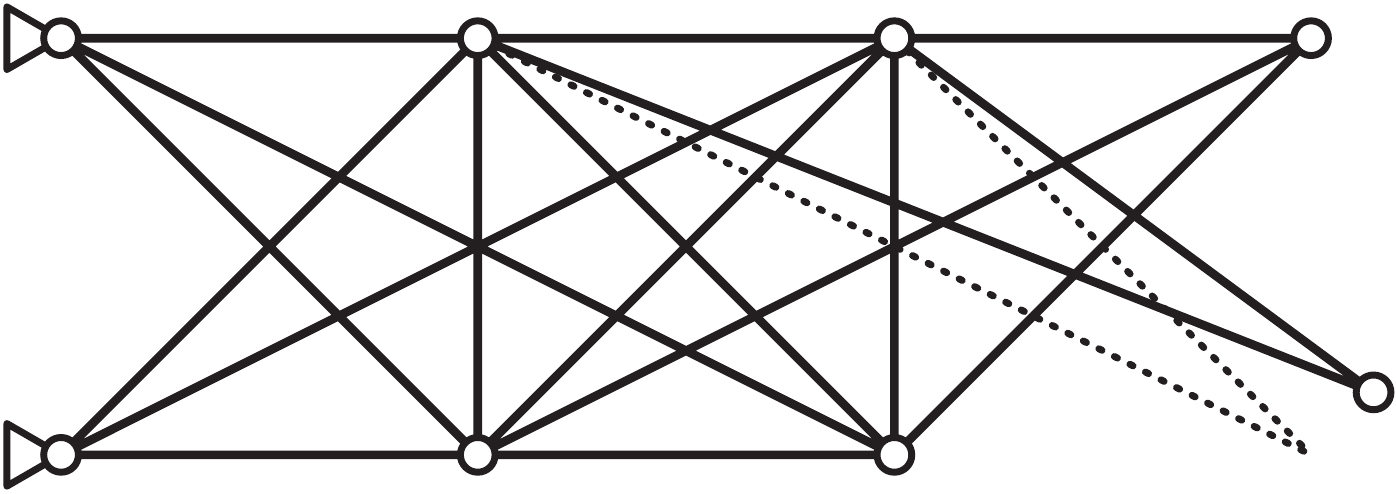}
    %\par
    %\includegraphics[scale=0.40]{ex2_enum3_9_alpha2.pdf}
    \caption{}
    \label{fig.ex2_enum3_9_alpha2}
  \end{subfigure}
  \par\bigskip
  \begin{subfigure}[b]{0.45\textwidth}
    \centering
    \includegraphics[scale=0.40]{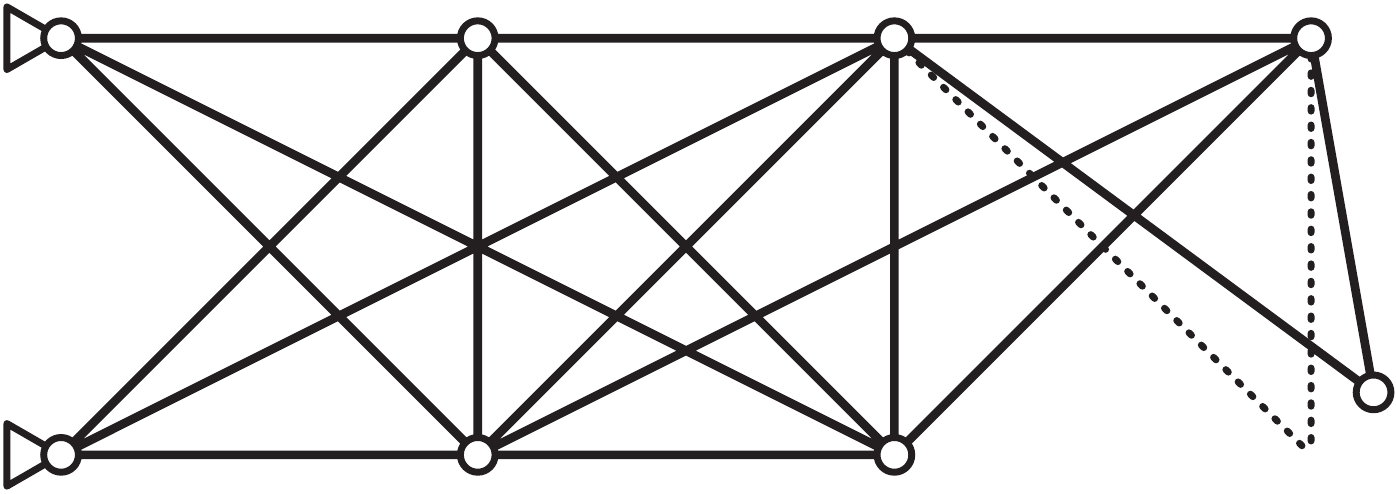}
    %\par
    %\includegraphics[scale=0.40]{ex2_enum3_19_alpha2.pdf}
    \caption{}
    \label{fig.ex2_enum3_19_alpha2}
  \end{subfigure}
  \hfill
  \begin{subfigure}[b]{0.45\textwidth}
    \centering
    \includegraphics[scale=0.40]{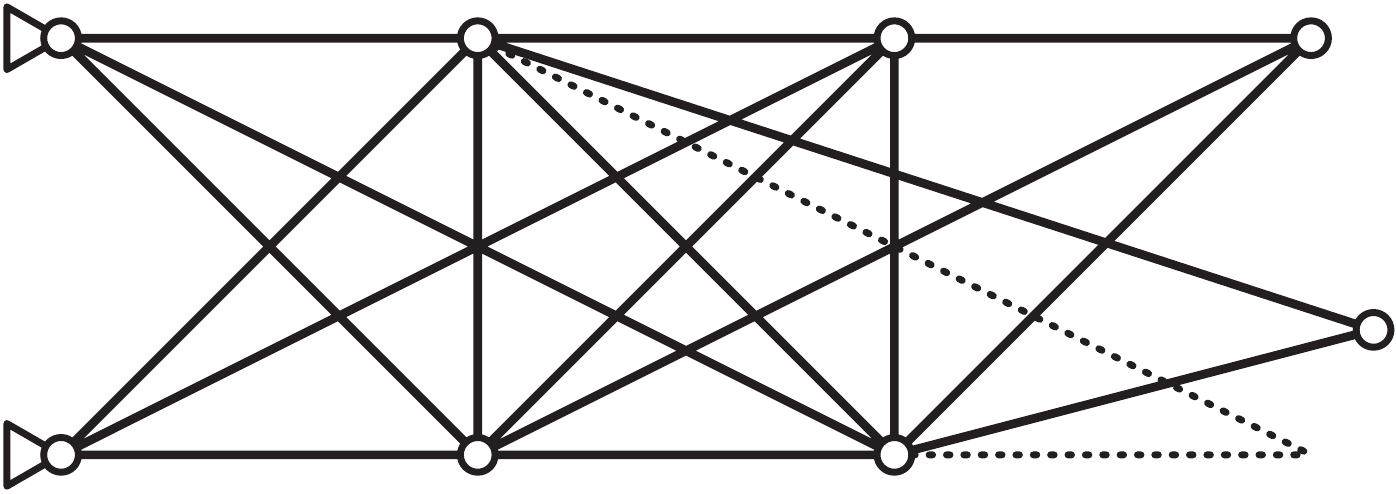}
    %\par
    %\includegraphics[scale=0.40]{ex2_enum9_15_alpha2.pdf}
    \caption{}
    \label{fig.ex2_enum9_12_alpha2}
  \end{subfigure}
  \par\bigskip
  \begin{subfigure}[b]{0.45\textwidth}
    \centering
    \includegraphics[scale=0.40]{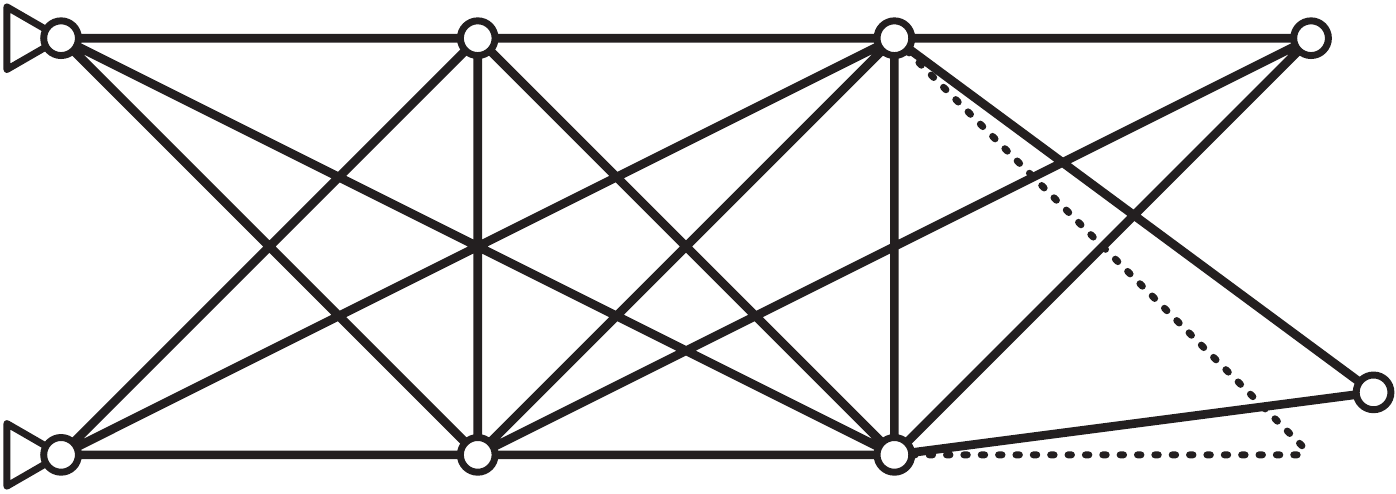}
    %\par
    %\includegraphics[scale=0.40]{ex2_enum9_19_alpha2.pdf}
    \caption{}
    \label{fig.ex2_enum9_17_alpha2}
  \end{subfigure}
  \hfill
  \begin{subfigure}[b]{0.45\textwidth}
    \centering
    \includegraphics[scale=0.40]{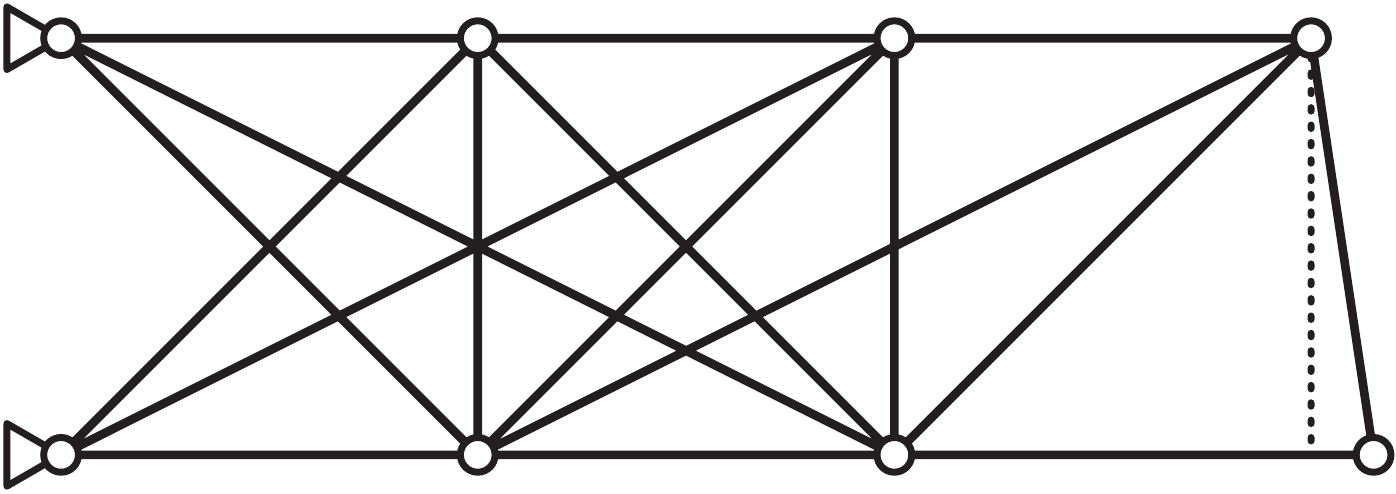}
    %\par
    %\includegraphics[scale=0.40]{ex2_enum15_19_alpha2.pdf}
    \caption{}
    \label{fig.ex2_enum12_17_alpha2}
  \end{subfigure}
  \par\bigskip
  \begin{subfigure}[b]{0.45\textwidth}
    \centering
    \includegraphics[scale=0.40]{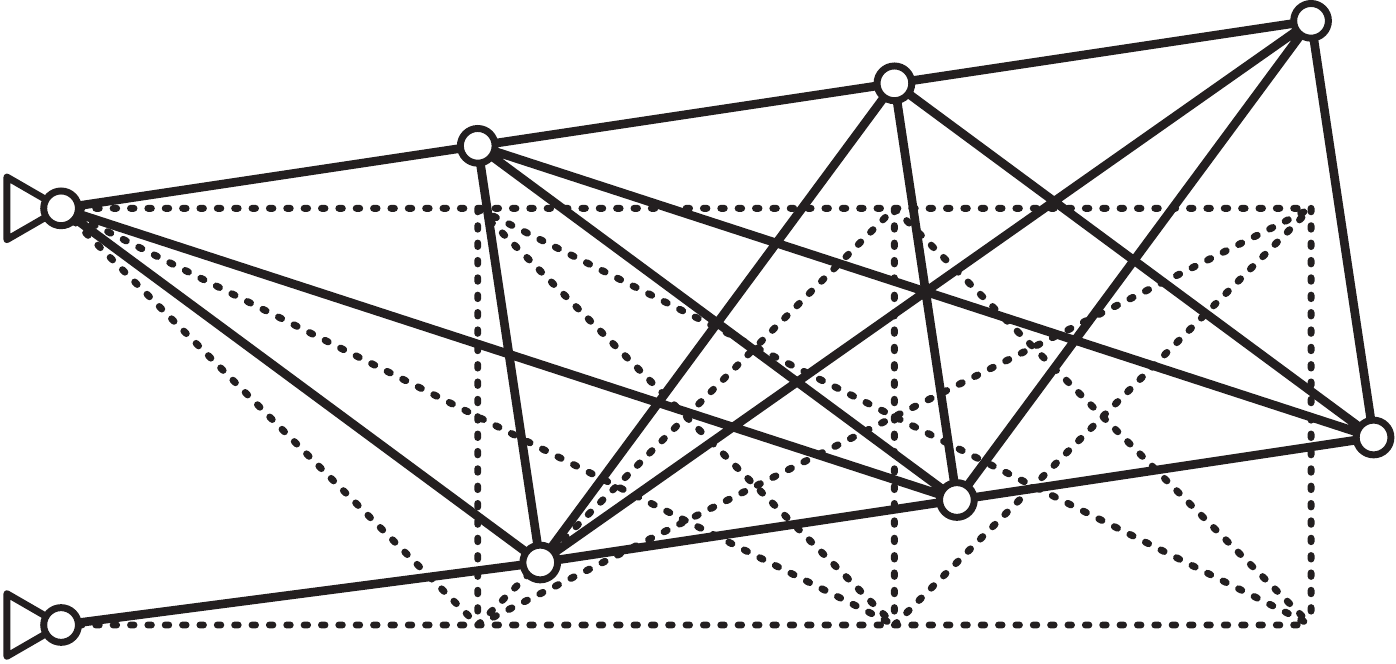}
    \caption{}
    \label{fig.ex2_enum13_18_alpha2}
  \end{subfigure}
  \hfill
  \begin{subfigure}[b]{0.45\textwidth}
  \end{subfigure}
  \caption{Some of the worst-case scenarios for the solution with 
  $\alpha=2$ in example (II). 
  The damage scenarios obtained by reflecting these scenarios across the 
  axis of symmetry are also the worst-scenarios. 
%  それぞれに対称なモードがあるので，重複度は 18 のはず（誤差により実際は 17）.
%  \subref{fig.ex2_enum13_18_alpha2} と対称なのは 
%  ``10 \& 16'' で $\lambda=3.408$. 
  }
  \label{fig.ex2_enum_alpha2}
\end{figure}

We next consider a different loading condition shown in 
\reffig{fig.m_19bar_2}, where proportionally increasing forces 
of $50\lambda\,\mathrm{kN}$ are applied horizontally at the two 
rightmost nodes. 
The optimal solution of the classical limit design problem (i.e., 
$\alpha=0$) is apparently the one shown in 
\reffig{fig.ex2_optimal_alpha_0}. 
Thus the optimal solution of the conventional optimization problem has 
no redundancy. 

When we set $\alpha=1$, \refalg{alg.SQP} finds the solution shown 
in \reffig{fig.ex2_optimal_alpha_1}. 
The algorithm terminates after solving $200$ QP problems and $1960$ MILP 
problems. 
At the final iteration, the stencil radius is 
$r=4.2 \times 10^{-3}\,\mathrm{mm^{2}}$, and 
the solution satisfies the termination condition at 
step~\ref{alg.SQP.subproblem}. 
The worst-case limit load factor of the obtained solution is $7.2812$, 
while that of the initial design is $5.7889$. 
At the obtained solution, the multiplicity of the worst-case scenarios 
is $9$. 
Among them, five scenarios are shown in \reffig{fig.ex2_enum_alpha1}; 
namely, since the obtained solution has symmetry, the damage scenarios 
obtained by reflecting the ones in Figures~\ref{fig.ex2_enum1_alpha1}, 
\subref{fig.ex2_enum3_alpha1}, \subref{fig.ex2_enum16_alpha1}, and 
\subref{fig.ex2_enum17_alpha1} across the axis of symmetry are also the 
worst-case scenarios. 
The collapse mode in \reffig{fig.ex2_enum9_alpha1} should be symmetric, 
if the obtained design is strictly symmetric. 
It is worth noting that we do not impose the constraints on symmetry of 
the design variables in the process of optimization. 
Asymmetry in \reffig{fig.ex2_enum9_alpha1} is due to numerical errors 
in symmetry. 

For $\alpha=2$, the solution obtained by \refalg{alg.SQP} is shown in 
\reffig{fig.ex2_optimal_alpha_2}, where $378$ QP problems and $4014$ 
MILP problems are solved. 
The algorithm terminates at step~\ref{alg.step.filter.outer.2}, which 
means that the objective value at any sample point is not better than 
that at the obtained solution. 
It is worth noting that the limit load factors of many damage scenarios 
coincide at the obtained solution. 
Indeed, the multiplicity of the worst-case scenarios is $18$. 
Among them, $9$ scenarios are collected in \reffig{fig.ex2_enum_alpha2}. 
Due to symmetry of the solution, the damage scenarios that are obtained 
by reflecting the ones in \reffig{fig.ex2_enum_alpha2} across the axis 
of symmetry are also the worst-case scenarios. 
The worst-case limit load factor of the obtained solution is $3.2773$, 
while that of the initial solution is $1.7889$. 

In all the examples presented above, all the candidate members of the 
ground structure present in the obtained solution. 
Since the proposed algorithm is based upon the SQP method, it in 
principle allows some members to vanish. 
Since the global optimality of the obtained solution is not guaranteed, 
it is possible that the obtained solution is only a local optimal 
solution that is not globally optimal and some members vanish at a 
global optimal solution. 
Or, a global optimal solution may truly have all the members in the 
ground structure. 
This issue remains to be studied.

\section{Summary and discussion}

This paper has defined a concept of redundancy optimization of 
structures. 
Roughly speaking, the redundancy optimization maximizes the structural 
functionality in the worst-case scenario when deficient structural 
components are unknown a priori. 
The notion of redundancy is related to robustness against uncertainty in 
the set of deficient components. 
In accordance with this relation, the proposed redundancy optimization 
formulation is naturally consistent with a robust optimization of 
structures, in which one attempts to optimize the worst-case value of 
the objective function under uncertainty in structural environment. 

A derivative-free optimization method has been proposed to 
solve the redundancy optimization problem. 
The method combines the SQP method and the finite-difference method with 
a varying difference increment. 
The numerical examples show that this algorithm can find a solution 
with multiple worst-case scenarios. 
Like a multiple eigenvalue in the eigenvalue optimization \citep{SLO94}, 
a multiple limit load factor may not be differentiable. 
Hence, it is rather surprising that the proposed SQP-based algorithm 
can find a solution with large multiplicity of limit load factors. 
Fundamental properties, such as continuity and smoothness, of the 
worst-case limit load factor as a function of the member cross-sectional 
areas are not investigated yet. 
Hence, the optimality condition of the redundancy optimization problem 
considered in this paper also remains to be studied. 

As a related design optimization problem, we may consider direct 
maximization of a quantitative measure of redundancy. 
For instance, the strong redundancy is defined by 
\begin{align*}
  \hat{\alpha}(\bi{x}; h^{\rr{c}}) 
  = \max \{ \alpha
  \mid
  h(\bi{s}) \le h^{\rr{c}} \ (\forall \bi{s} \in D(\bi{x}; \alpha))
  \} , 
\end{align*}
where $h^{\rr{c}}$ is the specified allowance of the structural 
performance \citep{KBh11}. 
Maximization of the strong redundancy is formulated as follows: 
\begin{subequations}\label{P.max.alpha.hat}
  \begin{alignat}{3}
    & \displaystyle 
    \Max
    &{\quad}& \displaystyle
    \hat{\alpha}(\bi{x}; h^{\rr{c}}) \\
    & \st && \displaystyle
    \bi{x} \in X . 
  \end{alignat}
\end{subequations}
Obviously, this problem is closely related to problem 
\eqref{P.redundancy.optimization}. 
Let $\bi{x}^{*}$ and $h^{*}$ denote the optimal solution and the optimal 
value of problem \eqref{P.redundancy.optimization}, respectively. 
If $h^{\rr{c}} \ge h^{*}$, then $\bi{x}^{*}$ satisfies 
$\hat{\alpha}(\bi{x}^{*}; h^{\rr{c}}) \ge \alpha$. 
Therefore, the optimal solution of problem \eqref{P.max.alpha.hat} can 
be explored by solving problem \eqref{P.redundancy.optimization} with 
varying the value of $\alpha$. 
Another possible formulation, that allows to handle several 
measures of structural performance, may be written as follows: 
\begin{subequations}\label{P.some.constraints}
  \begin{alignat}{3}
    & \displaystyle 
    \Min
    &{\quad}& \displaystyle
    c(\bi{x}) \\
    & \st && \displaystyle
    h_{j}(\bi{s}) \le h^{\rr{c}}_{j} 
    \ (\forall \bi{s} \in D(\bi{x}; \alpha)) , 
    \quad j=1,\dots,k , \\
    & && \displaystyle
    \bi{x} \in X . 
  \end{alignat}
\end{subequations}
Here, $c(\bi{x})$ is the cost function such as the structural volume. 
This optimization problem is consistent with the general methodology of 
robust optimization \citep{BtEgN09}. 
Also, problem~\eqref{P.some.constraints} with $\alpha=1$ is essentially 
similar to the fail-safe optimization problem of structures studied 
in \citet{SAH76} and \citet{NA82}. 

This paper has developed a generic framework for optimizing structures 
with guaranteeing the magnitude of redundancy. 
Other concepts of redundancy optimization based on 
different definitions of redundancy may be formulated. 
Although attention of this paper has been focused on truss structures, 
the presented concept can be applied to frame structures in a 
straightforward manner. 
Also, structural performance other than the limit load factor, such as 
the compliance and the violation of the displacement constraints, can be 
considered within the presented framework. 
The proposed algorithm might probably have many possibilities of 
improvements and extensions. 
Especially, improvements from the viewpoints of computational cost and 
convergence property may remain to be studied. 
Also, extension to nonlinear constraints could be explored. 
Moreover, in the presented numerical examples, all the solutions have 
all the members in the ground structure, although an SQP method 
generally allows some members to vanish if such a solution is optimal. 
It remains to be studied if the obtained solution is only a local 
optimal solution and some members vanish at a global optimal solution, 
or, a global optimal solution truly has all the members in the ground 
structure.

\section*{Acknowledgments}

This work is partially supported by 
the Support Program for Urban Studies from the Obayashi Foundation. 
%and 
%JSPS KAKENHI (C) 26420545. 

\end{document}